\documentclass[11pt,notitlepage,twoside]{article}
\usepackage{latexsym}
\usepackage{amsmath}
\usepackage{amsfonts}
\usepackage{amssymb}
\renewcommand{\a }{\alpha }

\renewcommand{\d}{\delta }

\newcommand{\D }{\Delta }

\newcommand{\e }{\varepsilon }

\renewcommand{\l }{\lambda }

\newcommand{\n }{\nabla }
\newcommand{\var }{\varphi }

\newcommand{\Sig }{\Sigma}

\renewcommand{\O }{\Omega }

\newcommand{\ov}{\overline}
\newcommand{\wtilde }{\widetilde}

\newcommand{\be}{\begin{equation}}
\newcommand{\ee}{\end{equation}}
\newenvironment{pf}{\noindent{\bf Proof.}\enspace}{
\hfill$\Box$\medskip}
\newenvironment{pfn}[1]{\noindent{\bf Proof of {#1}\enspace}}{
\hfill$\Box$\medskip}
\newcommand{\R}{\mathbb{R}}
 
\newcommand{\N}{\mathbb{N}}

\newtheorem{thm}{Theorem}[section]
\newtheorem{pro}[thm]{Proposition}
\newtheorem{lem}[thm]{Lemma}
\newtheorem{rem}[thm]{Remark}
\newtheorem{cor}[thm]{Corollary}

\numberwithin{equation}{section}

\author{Mohamed Ben Ayed$^a$\thanks{E-mail: \texttt{Mohamed.Benayed@fss.rnu.tn}.}, \, Khalil El Mehdi$^{b,c}$\thanks{E-mails:  \texttt{khalil@univ-nkc.mr}, \texttt{elmehdik@ictp.trieste.it}.} \& Mokhless Hammami$^a$ \thanks{Corresponding author. E-mail: \texttt{Mokhless.Hammami@fss.rnu.tn}. Fax: + 216-74 274437.}  \\
{\footnotesize
a: D{\'e}partement de Math{\'e}matiques, Facult{\'e} des Sciences de Sfax, Route
Soukra, Sfax, Tunisia}\\
{\footnotesize
b:  Facult\'e des Sciences et Techniques, Universit\'e de Nouakchott, Nouakchott, Mauritania}\\
{\footnotesize
 c: The Abdus Salam ICTP, Mathematics Section, Strada Costiera 11, 34014 Trieste, Italy.}
}

\title { \Large \textbf{Some Existence Results for a Paneitz Type Problem\\
Via the Theory of Critical Points at Infinity}}

\begin{document}

\date{ }
\maketitle
{\footnotesize
\noindent {\bf Abstract. } In this paper a fourth order equation  involving critical growth is considered under the Navier boundary condition:
 $\D ^2 u=K u^p$, $u>0$ in $\Omega$, $u=\D u=0$ on $\partial\Omega$,
where $K$ is a positive function, $\Omega$ is a bounded smooth
 domain in $\R^n$, $n\geq 5$ and  $p+1={2n/(n-4)}$ is the
critical Sobolev exponent. We give some topological conditions on
$K$ to ensure the existence of solution. Our methods involve the study of the critical points at infinity and their contribution to the topology of the level sets of the associated Euler Lagrange functional. \\
\noindent{\bf Mathematics Subject Classification (2000):} \quad 35J60, 35J65, 58E05.\\
\noindent{\bf Key words: } Critical
points at infinity, Critical Sobolev exponent, Lack of compactness.\\ 
\noindent{\bf R\'esum\'e.} Dans ce papier, nous consid\'erons une
\'equation d'ordre quatre ayant un accroissement critique avec
conditions de Navier au bord: $\D ^2 u=K u^p$, $u>0$ dans
$\Omega$, $u=\D u=0$ sur $\partial\Omega$, o\`u  $K$ est une
fonction strictement positive, $\Omega$ est un domaine born\'e
r\'egulier de $\R^n$, $n\geq 5$ et $p+1={2n/(n-4)}$ est l'exposant
critique de Sobolev. Nous donnons certaines conditions
topologiques sur $K$ pour assurer l'existence de solution. Notre
approche est bas\'ee sur l'\'etude des points critiques \`a
l'infini et de leur contribution \`a la topologie des ensembles de
niveau de la fonctionnelle d'Euler Lagrange associ\'ee.\\
\noindent{\bf Mots cl\'es: } Points critiques \`a l'infini,
Exposant critique de Sobolev, D\'efaut de compacit\'e. }
\section{Introduction  and Main Results}
\mbox{}
In this paper we prove some existence results for the following
nonlinear  problem under the Navier boundary condition
$$
(P)\qquad \left\{\begin{array}{cc}
\Delta ^2 u  =  K  u^p,\,  u > 0 &\mbox {in}\,  \Omega\\
\D u  =  u  =0 &\mbox {on}\,  \partial \Omega,
\end{array}
\right.
$$
where $\O$ is a bounded smooth domain of $\R^n$,
 $n\geq 5$, $p+1=\frac{2n}{n-4}$ is the critical exponent of the
embedding $H^2\cap H_0^1(\O)$ into $L^{p+1}(\O)$ and $K$ is a $C^3$-positive function in $\ov{\O}$.

This type of equation naturally arises from the study of conformal
geometry. A well known example is the problem of prescribing the
Paneitz curvature : given a function $K$ defined in compact
Riemannian manifold $(M,g)$ of dimension $n\geq 5$, we ask whether
there exists a metric $\tilde{g}$ conformal to $g$ such that $K$
is the Paneitz curvature of the new metric $\tilde{g}$ (for
details one can see \cite{BE1}, \cite{BE2}, \cite{C},  \cite{DHL},
\cite{DMO1}, \cite{DMO2}, \cite{F} and the references
therein).

We observe that one of the main features of problem $(P)$ is the lack of compactness, that is, the Euler Lagrange functional $J$ associated to $(P)$ does not satisfy the Palais-Smale condition. This means that there exist noncompact sequences along which the functional is bounded and its gradient goes to zero. Such a fact follows from the noncompactness of the embedding of   $H^2\cap H_0^1(\O)$ into $L^{p+1}(\O)$. However, it is easy to see that a necessary condition for solving the problem $(P)$ is that $K$ has to be positive somewhere. Moreover, it turns out that there is at least another obstruction to solve the problem $(P)$, based on Kazdan-Warner type conditions, see \cite{DHL}. Hence it is not expectable to solve problem $(P)$ for all the functions $K$, thus a natural question arises: under which conditions on $K$, $(P)$ has a solution. Our aim in this paper is to give sufficient conditions on $K$ such that $(P)$ possesses a solution.

In the last years, serval researches have been developed on the existence of solutions of fourth order elliptic equations with critical exponent on a domain of $\R^n$, see \cite{BH}, \cite{BGP}, \cite{EFJ}, \cite{EO}, \cite{GGS}, \cite{HV}, \cite{Lin}, \cite{NSJ}, \cite{PV}, \cite{PS}, \cite{V1} and \cite{V2}. However, at the authors' knowledge, problem $(P)$ has been considered for $K\equiv 1$ only.

As we mentioned before, $(P)$ is delicate from a variational viewpoint because of the failure of the Palais-Smale condition, more precisely because of the existence of critical points at infinity, that is orbits of the gradient flow of $J$ along which $J$ is bounded, its gradient goes to zero, and which do not converge \cite{B1}.
In this article, we give a contribution in the same direction as in the papers \cite{AB}, \cite{B2}, \cite{BCH} concerning  the problem of prescribing the scalar curvature on closed manifolds. Precisely, we extend some topological and dynamical methods of the {\it Theory of critical points at infinity} (see \cite{B1}) to the framework of such higher order equations. To do such an extension, we perform a careful expansion of $J$, and its gradient near a neighborhood of higly concentrated functions. Then, we construct a special pseudogradient for the associated variational problem for which the Palais-Smale condition is satisfied along the decreasing flow lines far from a finite number of such ``singularities''. As a by product of the construction of our pseudogradient, we are able to characterize the critical points at infinity of our problem. Such a fine analysis  of these critical points at infinity, which has its own interest, is highly nontrivial and plays a crucial role in the derivation of  existence results. In our proofs, the main idea is  to take advantage of the precise computation of the contribution of these critical points at infinity to the topology of the level sets of $J$; the main argument being that, under our conditions on $K$, there remains some difference of topology which is not due to the critical points at infinity and therefore the existence of a critical point of $J$.

Our proofs go along the methods of Aubin-Bahri \cite{AB}, Bahri \cite{B2} and  Ben Ayed-Chtioui-Hammami \cite{BCH}. However, in our case the presence of the boundary makes the analysis more involved: it turns out that the interaction of ``bubbles'' and the boundary creates a phenomenon of new type which is not present in the closed manifolds'case. In addition, we have to prove the positivity of the critical point obtained by our process.  It is known that in the framework of  higher order equations such a proof is quite difficult in general (see \cite{DMO2} for example), and the way we handle it here is very simple compared with the literature, see Proposition \ref{p:40} below.

In order to state our main results, we need to introduce some notation and the assumptions that we are using in our results.
We denote by $G$ the Green's function and by $H$ its regular part, that is
 for each $x\in \O$,
$$
 \left\{\begin{array}{ccccc}
G(x,y) & = & \mid x-y\mid^{-(n-4)}-H(x,y)  &  \mbox { in }\,  \O,  \\
\D ^2 H(x,.)& = & 0  & \mbox { in }\, \O , \\
\D G(x,.) & = & G(x,.) =0 & \mbox{ on } \partial \O.
\end{array}
\right.
$$
Now, we state our assumptions.\\
 $(A_0)$ \quad Assume that, for each $x\in \partial\O$
$$\frac{\partial K(x)}{\partial \nu}<0,
$$
where $\nu$ is the outward normal to $\O$.\\
 $(A_1)$\quad We assume that  $K$ has only nondegenerate
critical points $y_0$, $y_1$,..., $y_s$ such that
$$
K(y_0)\geq K(y_1)\geq ...\geq K(y_l)>K(y_{l+1})\geq ...\geq
K(y_s).
$$
 $(A_2)$\quad We assume that
 $$
-\frac{\D K(y_i)}{60K(y_i)}+H(y_i,y_i) > 0 \,\,
\mbox{for}\, i\leq l \,\, \mbox{and} \,\, -\frac{\D
K(y_i)}{60K(y_i)}+H(y_i,y_i) < 0 \,\, \mbox{for}\, i > l\,
(\mbox{if } n=6),
 $$ 
\hskip 1cm and
$$
-\D K(y_i)>0 \, \mbox{ for } i\leq l\quad \mbox{ and }
-\D K(y_i)<0 \, \mbox{ for } i> l\quad (\mbox{ if }n\geq 7).
$$
 $(A'_2)$\quad  We assume that
$$
 -\frac{1}{60}\frac{\D
K(y_i)}{K(y_i)}+H(y_i,y_i) < 0 \,\, \mbox{for}\,  i > l\,\,
(\mbox{if } n=6) \, \mbox{ and }\, -\D K(y_i)<0 \, \mbox{ for } i>
l\,  (\mbox{if }n\geq 7).
$$
In addition, for every $i\in\{1,...,l\}$ such that
$$
-\frac{1}{60}\frac{\D K(y_i)}{K(y_i)}+H(y_i,y_i) \leq 0 \,\,
(\mbox{if } n=6) \,  \mbox{ and }\,  -\D K(y_i)\leq 0 \quad
(\mbox{if }n\geq 7),
$$
we assume that
 $n-m+3 \leq \mbox{index}(K,y_{i})\leq n-2$,
where index$(K,y_{i})$ is the Morse index of $K$ at $y_{i}$
and $m$ is an integer defined in assumption $(A_3)$.\\
Now, let $Z_K$ be a pseudogradient of $K$ of Morse-Smale
type (that is, the intersections of the stable and unstable
manifolds of the critical points of $K$ are transverse). Set
 $$X=\ov{\bigcup_{0\leq i\leq l}W_s(y_i)},$$
 where $W_s(y)$ is the stable manifold of $y$ for $Z_K$.\\
 $(A_3)$\quad We assume that $X$ is not contractible and denote by $m$ the
 dimension of the first nontrivial reduced homological group of $X$.\\
  $(A_4)$ \quad We assume that there exists a positive constant $\ov{c}<K(y_l)$ such that
  $X$ is contractible in $K^{\ov{c}}=\{x\in \O/K(x)\geq \ov{c}\}$.

Now we are able to state our first results
\begin{thm}\label{t:11}
 Let $n\geq 6$. Under the assumptions $(A_0)$, $(A_1)$, $(A_2)$, $(A_3)$ and $(A_4)$,
there exists a constant $c_0$ independent of $K$ such that if
$K(y_0)/\ov{c}\leq 1+c_0$, then (P) has a solution.
  \end{thm}
\begin{cor}\label{c:12}
The solution obtained in Theorem \ref{t:11} has an augmented Morse index $\geq m$.
\end{cor}
\begin{thm}\label{t:13}
Let $n\geq 7$. Under the assumptions $(A_0)$, $(A_1)$, $(A'_2)$,
$(A_3)$ and $(A_4)$, there exists a constant $c_0$ independent of
$K$ such that if $K(y_0)/\ov{c}\leq 1+c_0$, then (P) has a
solution.
  \end{thm}
\begin{rem}
{\bf i).} The assumption $K(y_0)/\bar{c}\leq 1+c_0$ allows basically to perform a single-bubble analysis.\\
{\bf ii).} To see how to construct an example of a function $K$ satisfying our assumptions, we refer the interested reader to \cite{AB2}.
\end{rem}
Next, we state another kind of existence results for problem $(P)$
based on a topological invariant introduced by A. Bahri in
\cite{B2}. In order to give our results in this direction, we need
to fix some notation and state our assumptions.\\
We denote by $W_s(y)$ and $W_u(y)$ the stable and unstable
manifolds of $y$ for $Z_K$. \\
$(A_5)$\quad We assume that $K$ has only nondegenerate critical points $y_i$
satisfying $\D K(y_i)\ne 0$ and \\$W_s(y_i)\cap
W_u(y_j)=\emptyset$ for any $i$ such that
$-\D K(y_i)>0$ and for any $j$ such that \\$-\D K(y_j)<0$.\\
For $k\in \{1,...,n-1\}$, we define $X$ as
$$
X=\ov{W_s(y_{i_0})},
$$
where $y_{i_0}$ satisfies
$$
K(y_{i_0})= \mbox{max }\{K(y_i)/\mbox{index }(K,y_i)=n-k, \quad
-\D K(y_i) >0 \}.
$$
$(A_6)$\quad  We assume that $X$ is without boundary.\\
We observe that assumption $(A_0)$ implies that $X$ does not
intersect the boundary $\partial\O$ and therefore it is a compact set of $\O$.\\
Now, we denote by $y_0$ the absolute maximum of $K$. Let us define the
set $C_{y_0}(X)$ as
$$C_{y_0}(X)=\{\a\d_{y_0}+(1-\a)\d_x/\a\in [0,1],\, x\in X\},
$$
where $\d_x$ denotes the Dirac mass at $x$.\\
 For $\l$ large enough,
we introduce a map $f_{\l}: C_{y_0}(X)\rightarrow \Sig^+:=\{ u \in
H^2\cap H_0^1/u>0, ||u||_2=1\} $
$$\a\d_{y_0}+(1-\a)\d_x\rightarrowtail
\frac{(\a/K(y_0)^{(n-4)/8})P\d_{(y_0,\l)}+((1-\a)/K(x)^{(n-4)/8})
P\d_{(x,\l)}}{\mid\mid (\a/K(y_0)^{(n-4)/8})P\d_{(y_0,\l)}
+((1-\a)/K(x)^{(n-4)/8})P\d_{(x,\l)}\mid\mid_2},
$$
where $||u||^2_2 = \int_\O |\D u|^2$.\\
 Then $C_{y_0}(X)$ and $f_{\l}(C_{y_0}(X))$ are manifolds in
 dimension $k+1$, that is, their singularities arise in dimension
$k-1$ and lower, see \cite{B2}. The codimension of
$W_s(y_0,y_{i_0})_\infty$ is equal to $k+1$, then we can define
the intersection number (modulo 2) of $f_{\l}(C_{y_0}(X))$ with
$W_s(y_0,y_{i_0})_\infty$
 $$\mu(y_{i_0})=f_{\l}(C_{y_0}(X)).W_s(y_0,y_{i_0})_\infty,$$
where $W_s(y_0,y_{i_0})_\infty$ is the stable manifold of the critical point at infinity 
$(y_0,y_{i_0})_\infty$ for a decreasing pseudogradient for $J$
which is transverse to $f_{\l}(C_{y_0}(X))$. Such a number is
well defined  see \cite{B2},\cite{M}. Observe that $C_{y_0}(X)$ and
 $f_{\l}(C_{y_0}(X))$ are contractible while $X$ is not
 contractible.\\
$(A_7)$\quad Assume that $2/K(y_0)^{(n-4)/4}<1/K(y)^{(n-4)/4}$ for
each critical point $y$ of Morse index $n-(k+1)$ and satisfies
$-\D K(y)>0$.\\
 We then have the following result:
\begin{thm}\label{t:14}
Let $n\geq 7$. Under assumptions $(A_0)$, $(A_5)$, $(A_6)$ and $(A_7)$, if
$\mu(y_{i_0})=0$ then (P) has a solution of an augmented Morse index less than $k+1$.
\end{thm}
Now, we give a more general statement than Theorem \ref{t:14}.
For this purpose, we define $X$ as
$$
X=\ov{\cup_{y\in B}W_s(y)},
$$
where $B=\{y\in \O/\n K(y)=0,\, -\D K(y)>0\}$. We denote by
$k$ the dimension of $X$ and by $B_k=\{y\in B/index(K,y)=n-k\}$.\\
For $y_i \in B_k$, we define, for $\l$ large enough, the intersection number (modulo $2$)
 $$\mu(y_i)=f_{\l}(C_{y_0}(X)).W_s(y_0,y_i)_\infty.$$
By the above arguments, this number is well defined, see \cite{M}.\\
Then, we have:
\begin{thm}\label{t:15}
Let $n\geq 7$. Under assumptions $(A_0)$, $(A_5)$ and $(A_6)$, if
$\mu (y_i)=0$ for each $y_i\in B_k$, then (P) has a solution of
an augmented Morse index less than $k+1$.
\end{thm}

The organization of the paper is the following.
In section 2, we set up the variational structure and recall some
preliminaries. In section 3, we give an expansion of the Euler
functional associated to $(P)$ and its gradient near potential
critical points at infinity. In section 4, we provide the proof of
Theorem \ref{t:11} and its corollary. In section 5, we prove
Theorem \ref{t:13}, while section 6 is devoted to the proof of
Theorems \ref{t:14} and \ref{t:15}.
\section{Preliminaries}
\mbox{}
In this section, we set up the variational structure and its mean features.\\
 Problem $(P)$ has a variational structure. The related functional
is
 $$ J(u) =\biggl(\int_\Omega\ K \mid
u\mid^{\frac{2n}{n-4}}\biggr)^{-\frac{n-4}{n}}$$
 defined on
$$
\Sigma =\{ u\in H^2\cap H_0^1(\O) / \mid\mid u\mid\mid_{H^2\cap
H_0^1(\O)}^2:=\mid\mid u\mid\mid_2^2:= \int_\O\mid\D u\mid ^2
=1\}.
$$
The positive critical points of $J$ are solutions of ($P$), up to
a
multiplicative constant.\\
 Due to the non-compactness of the
embedding $H^2\cap H_0^1(\O)$ into $L^{p+1}(\O)$, the functional
$J$ does not satisfy the Palais-Smale condition. An important result of
Struwe \cite{S} (see also \cite{L} and \cite{BrC}) describes the behavior of
such sequences associated to second order equations of the type
\begin{eqnarray}
-\D u =\mid u^{\frac{n+2}{n-2}},\quad u>0\quad \mbox{in}\quad \O;
\qquad u=0\quad\mbox{on}\quad \partial\O.
\end{eqnarray}
In \cite{GGS}, Gazzola, Grunau and Squassina proved the analogue  of this result for  problem $(P)$.
To
describe the sequences failing the Palais-Smale condition, we
need to introduce some notation.\\
 For $a\in\Omega $  and $ \lambda > 0 $, let
\begin{eqnarray}\label{d}
\delta_{(a,\lambda )}(x)= c_n \biggl (\frac{\lambda} {1+\lambda
^2\mid x-a\mid^2}\biggr ) ^{\frac{n-4}{ 2}},
\end{eqnarray}
where $c_n$ is a positive constant chosen so that
$\delta_{(a,\lambda)}$ is the family of solutions of the following
problem (see \cite{Lin})
\begin{eqnarray}
\D ^2 u=\mid u^{\frac{n+4}{n-4}}, \quad u>0\quad \mbox{ in }
\R^n.
\end{eqnarray}
 For $f\in H^2(\O)$, we define the projection $P$ by
\begin{eqnarray}
u=Pf\Longleftrightarrow\D ^2u=\D^2f \mbox{ in } \O,\quad u=\D u=0
\mbox{ on } \partial \O.
\end{eqnarray}
We have the following proposition which is extracted from
\cite{BH}.
\begin{pro}\label{13}\cite{BH}
Let $a\in \O$, $\l>0$ and
$\varphi_{(a,\l)}=\d_{(a,\l)}-P\d_{(a,\l)}$. We have
$$
(a)\quad 0\leq \varphi_{(a,\l)}\leq\d_{(a,\l)}, \qquad (b) \quad
\varphi_{(a,\l)}=c_n\frac{H(a,.)}{\l^{\frac{n-4}{2}}}+f_{(a,\l)}$$
where
 $c_n$ is defined in \eqref{d} and $f_{(a,\l)}$ satisfies
$$f_{(a,\l)}=O\left(\frac{1}{\l^{\frac{n}{2}}d^{n-2}}\right),\quad
\l\frac{\partial f_{(a,\l)}}{\partial\l}=O\biggl(\frac{1}{
\l^{\frac{n}{2}}d^{n-2}}\biggr),\quad \frac{1}{\l}\frac{\partial
f_{(a,\l)}}{\partial a}=O\biggl(\frac{1}{
\l^{\frac{n+2}{2}}d^{n-1}}\biggr)$$ where $d$ is the distance
$d(a,\partial \O)$.
$$\mid\varphi_{(a,\l)}\mid_{L^{\frac{2n}{n-4}}}=O\bigl(\frac{1}{(\l
d)^{\frac{n-4}{2}}}\bigr), \quad
\mid\l\frac{\partial\varphi_{(a,\l)}}{\partial\l}\mid_{L^{\frac{2n}{n-4}}}=O\bigl(
\frac{1}{(\l d)^{\frac{n-4}{2}}}\bigr), \leqno{(c)}$$
 $$\mid\mid\varphi_{(a,\l)}\mid\mid_{2}=O\bigl(\frac{1}{(\l
d)^{\frac{n-4}{2}}}\bigr), \quad
 \mid\frac{1}{\l}\frac{\partial\varphi_{(a,\l)}}{\partial a}
\mid_{L^{\frac{2n}{n-4}}}=O\bigl(\frac{1}{(\l
d)^{\frac{n-2}{2}}}\bigr).$$
\end{pro}

We now introduce the set of potential critical points at infinity.\\
For any $\e>0$ and $p\in \N^*$, let $V(p,\e)$ be the subset of
$\Sig$ of the following functions: $u\in \Sig$ such that there is
$(a_1,...,a_p)\in \O^p$, $(\l_1,...,\l_p)\in (\e^{-1},+\infty)^p$
and $(\a_1,...,\a_p)\in (0,+\infty)^p$ such that
$$\bigg |\bigg |u-\sum_{i=1}^p\a_i
P\d_{(a_i,\l_i)}\bigg |\bigg |_{2}<\e,\,\,
\l_id(a_i,\partial\O)>\e^{-1},\,\, \bigg |
\frac{\a_i^{8/(n-4)}K(a_i)}{\a_j^{8/(n-4)}K(a_j)}-1\bigg | <\e,\,\,
\e_{ij}<\e \mbox{ for } i\ne j $$ where
\begin{eqnarray}
\e_{ij}=\biggl(\frac{\l_i}{\l_j}+\frac{\l_j}{\l_i}+\l_i\l_j\mid
a_i-a_j\mid^2\biggr)^{-\frac{n-4}{2}}.
\end{eqnarray}
The failure of the Palais-Smale condition can be described going along the ideas developed in \cite{BrC}, \cite{L}, \cite{S}. Namely, we have:
\begin{pro}\label{p:21} \cite{GGS} 
Assume that $J$ has no critical point in $\Sig^+$. Let ($u_k$)
$\in\Sigma ^+$ be a sequence such that $(\partial J(u_k))$ tends
to zero and $(J(u_k))$ is bounded. Then, after possibly having
extracted a subsequence, there exist $p\in N^*$ and a sequence
$(\e_k)$, $\e_k$ tends to zero, such that $u_k\in V(p,\e_k)$.
\end{pro}

Now, we consider the following minimization problem
for a function $u\in V(p,\e)$ with $\e$ small
\begin{eqnarray}\label{e:26}
\min \{\mid\mid
u-\sum_{i=1}^p\a_iP\d_{(a_i,\l_i)}\mid\mid_{2},\quad
\a_i>0,\quad\l_i>0,\quad a_i\in \O\}.
\end{eqnarray}
We then have the following proposition whose proof is similar, up
to minor modifications, to the corresponding statement for the
Laplacian operator in \cite{BC}. This proposition defines a
parametrization of the set $V(p,\e )$.
\begin{pro}\label{p:22}
For any $p\in \N^*$, there exists $\e_p>0$ such that, if $\e<\e_p$
and  $u\in V(p,\e)$, the minimization problem \eqref{e:26} has a
unique solution $(\a,a,\l)$ (up to permutation). In particular, we
can write $u\in V(p,\e )$ as follows
$$
u=\sum_{i=1}^p\a_iP\d_{(a_i,\l_i)} + v,
$$
where $(\a_1,...,\a_p,a_1,...,a_p,\l_1,...,\l_p)$ is the solution
of \eqref{e:26} and $v\in  H^2(\O)\cap H^1_0(\O)$ such that
\begin{eqnarray*}
(V_0)\qquad (v,P\d_{(a_i,\l_i)})_2=(v,\partial P\d_{(a_i,
\l_i)}/\partial \l_i)_2=0,\,\, (v,\partial P\d_{(a_i,
\l_i)}/\partial a_i)_2=0\,\, \mbox{for }i=1,...,p,
\end{eqnarray*}
where $(u,w)_2=\int_{\O}\D u\D w$.
\end{pro}

\section{  Expansion of the Functional and its Gradient }

 In this section, we will give a useful expansion of the functional $J$
and its gradient in the potential set $V(p,\e)$ for $n\geq 6$. In the sequel, for
the sake of simplicity, we will write $\d_i$ instead of
$\d_{(a_i,\l_i)}$. We start by the expansion of $J$.
 \begin{pro}\label{p:31}
There exists $\e_0>0$ such that for any
$u=\sum_{i=1}^p\alpha_iP\delta_i+v\in V(p,\e)$, $\e<\e_0$, $v$
satisfying $(V_0)$, we have
\begin{align*}
J(u)= & \frac{S_n^{4/n}\sum_{i=1}^p\a_i ^2}{(\sum_{i=1}^p\a_i
^{\frac{2n}{n-4}}
K(a_i))^{\frac{n-4}{n}}}\biggl[1+\frac{1}{S_n\sum_{i =1}^p
K(a_i)^{\frac{4-n}{4}}}\biggl(-\frac{n-4}{n}c_3\sum_{i=1}^p
\frac{\D K(a_i)}{K(a_i)^{n/4}\l_i ^2}\\
 & +c_2 \sum_{i=1}^p\frac{H(a_i,a_i)}{K(a_i)^{(n-4)/4}\l_i
^{n-4}}-\frac{c_2}{(K(a_i)K(a_j))^{(n-4)/8}}\sum_{i\ne j}\biggl(
\e_{ij}-\frac{H(a_i,a_j)}{(\l_i\l_j)^{(n-4)/2}} \biggr)
\biggr)\\
 & -f(v)+\frac{1}{\sum_{i=1}^p\a_i ^2S_n}Q(v,v)
+o\biggl(\sum \frac{1}{\l_k^2}+\frac{1}{(\l_kd_k)^{n-4}}+
\sum_{i\ne j} \e_{ij}+\mid\mid v\mid\mid_{2}^2\biggr)\biggr]
\end{align*}
where
 \begin{align*}
  Q(v,v)& =\mid\mid v\mid\mid_{ 2}^2-\frac{n+4}{n-4}
\sum_{i=1}^p\int_\O P\d_i ^{\frac{8}{n-4}}v^2\,  , \\
f(v)& =\frac{2}{\sum_{i=1}^p\a_i ^{2n/(n-4)} K(a_i)S_n}\int_\O K
\bigl(\sum_{i=1}^p\a_i P\d_i \bigr)^{\frac{n+4}{n-4}}v\quad ,
\end{align*}
$$S_n=\int_{\R^n}\frac{c_n^\frac{2n}{n-4} dy}{(1+ \mid
y\mid^2)^{n}},\quad
c_2=\int_{\R^n}\frac{c_n^{\frac{2n}{n-4}}}{(1+|y|^2)^{\frac{n+4}{2}}}dy\,
,\quad
c_3=\frac{c_n^{\frac{2n}{n-4}}}{2n}\int_{\R^n}\frac{|y|^2}{(1+|y|^2)^{n}}dy,$$
and $c_n$ is defined in \eqref{d}. Observe that if $n=6$ we have
$c_2=20c_3$.
\end{pro}
\begin{pf}
On  one hand, Proposition \ref{13} implies
\begin{align}
 \mid\mid P\delta\mid\mid_2^2&=
S_n-c_2\frac{H(a,a)}{\l ^{n-4}}+O\biggl(\frac{1}{(\l d)^{n-2}}\biggr),\label{o:33}\\
\int_\O K P\d ^{\frac{2n}{n-4}}&=K(a)S_n+c_3\frac{\D K(a)}{\l
^2}-\frac{2n}{n-4}c_2K(a)\frac{H(a,a)}{\l^{n-4}}+
O\biggl(\frac{1}{\l^3}+\frac{1}{(\l d)^{n-2}}\biggr)\label{o:34}.
\end{align}
On the other hand, a computation similar to the one performed in \cite{B1} shows that, for $i\ne j$, we have
\begin{eqnarray}\label{o:31}
 \int_{\R^n}\d_i ^{\frac{n+4}{n-4}}
\d_j=c_2\e_{ij}+ O(\e_{ij}^{\frac{n-2}{n-4}})\quad, \quad
\int_{\R^n}(\d_i\d_j) ^{\frac{n}{n-4}} =
O(\e_{ij}^{\frac{n}{n-4}}log(\e_{ij}^{-1})).
\end{eqnarray}
Thus, we derive that
\begin{eqnarray}\label{o:35}
\bigl(P\delta_i,P\delta_j\bigr)_2=c_2 \biggl(
\e_{ij}-\frac{H(a_i,a_j)}{(\l_i\l_j)^{(n-4)/2}}\biggr)+
O\biggl(\e_{ij}^{\frac{n-2}{n-4}} +\sum_{k=i,j}
\frac{1}{(\l_kd_k)^{n-2}}\biggr),
\end{eqnarray}
\begin{eqnarray}\label{o:36}
\int_\O K P\d_i^{\frac{n+4}{n-4}}P\d_j  = K(a_i) \bigl(P\d_i,P
\d_j\bigr)_2+ o\biggl(\sum
\frac{1}{\l_k^2}+\frac{1}{(\l_kd_k)^{n-4}}+ \e_{ij}\biggr)
\end{eqnarray}
and
\begin{eqnarray}\label{o:o}
\int K(\sum_{i=1}^p\a_iP\d_i)^{\frac{8}{n-4}}v^2=\sum_{i=1}^p
\a_i ^{\frac{8}{n-4}}K(a_i)\int P\d_i
^{\frac{8}{n-4}}v^2+o(\mid\mid v\mid\mid_{2}^2).
\end{eqnarray}
Combining \eqref{o:33},..., \eqref{o:o} and the fact that $\a_i
^{\frac{8}{n-4}}K(a_i)/(\a_j^{ \frac{8}{n-4}}K(a_j))=1+o(1)$, our
result follows.
\end{pf}\\
Now, let us recall that the quadratic form $Q(v,v)$  defined in
Proposition \ref{p:31} is positive definite (see \cite{BE1}). Thus
we have the following proposition which deals with the $v$-part of
$u$.
\begin{pro}\label{p:32}(see  \cite{BE1})
There exists a $C^1$-map which, to each $(\a, a, \l)$ satisfying\\
$\sum_{i=1}^p\a_iP\d_{(a_i,\l_i)}\in V(p,\e)$, with $\e$ small
enough, associates $\ov{v}=\ov{v}(\a,a,\l)$ satisfying $(V_0)$
such that $\ov{v}$ is unique, minimizing
$J(\sum_{i=1}^p\a_iP\d_{(a_i,\l_i)}+v)$ with respect to $v$
satisfying $(V_0)$, and we have the following estimate
\begin{align*}
\mid\mid\ov{v}\mid\mid_{2} \leq c\mid f\mid &
=O\biggl(\sum_{i=1}^p \frac{\mid\n K(a_i)\mid}{\l_i}+\frac{1}{\l_i
^2}\biggr) +(\mbox{if }n<12) O
\biggl(\sum{\e_{ij}(log\e_{ij}^{-1})^{\frac{n-4}{n}}}
+\frac{1}{(\l_id_i)^{n-4}}\biggr)\\
 & +(\mbox{if } n\geq
12)O\biggl(\sum \e_{ij}^{\frac{n+4}{2(n-4)}}
(log\e_{ij}^{-1})^{\frac{n+4}{2n}}+
\frac{(log\l_id_i)^{\frac{n+4}{2n}}}{(\l_id_i)^{\frac{n+4}{2}}}
\biggr).
\end{align*}
\end{pro}
Now regarding the  gradient of $ J$ which we will denote by $\partial J$, we have the following expansions
\begin{pro}\label{p:33}
 For $u=\sum_{i=1}^p \alpha_i
P\delta_i \in V(p,\varepsilon)$, we have the following expansion
\begin{align*}
\biggl( \partial J(u), \l_i\frac{\partial P\d_i}{\partial
\l_i}\biggr)_2= & 2J(u)\biggl[\frac{n-4}{n}c_3\a_i\frac{\D
K(a_i)}{K(a_i) \l_i ^2}-\frac{n-4}{2}c_2\a_i
\frac{H(a_i,a_i)}{\l_i ^{n-4}}\bigl(1+o(1)\bigr)\\
 & -c_2\sum_{j\ne i} \a_j\biggl(\l_i\frac{\partial\e_{ij}}{\partial
\l_i}+\frac{n-4}{2}\frac{H(a_i,a_j)}{(\l_i\l_j)^{(n-4)/2}}\biggr)\bigl(1+o(1)\bigr)
\biggr]\\
 & + o\biggl(\sum
\frac{1}{\l_k^2}+\frac{1}{(\l_kd_k)^{n-3}}+\sum_{k\ne r}
 \e_{kr}^{\frac{n-3}{n-4}}\biggr).
\end{align*}
\end{pro}
\begin{pf} We have
\begin{align}\label{e:31}
\biggl(\partial J(u),\lambda_i\frac{\partial
P\delta_i}{\partial\lambda_i} \biggr)_2= &
2J(u)\biggl[\sum\alpha_j\biggl(P\delta_j,\lambda_i\frac{\partial
P\delta_i}{\partial\lambda_i} \biggr)_2\\
 & -J(u)^{\frac{n}{n-4}}\int K
\bigl(\sum\alpha_jP\delta_j\bigr)^{\frac{n+4}{n-4}}\lambda_i\frac{\partial
P\delta_i}{\partial\lambda_i}\biggr].\notag
\end{align}
Observe that
\begin{align}\label{e:32}
\bigl(\sum&\alpha_j P\delta_j\bigr)^{\frac{n+4}{n-4}}  =
\sum\bigl(\alpha_jP\delta_j\bigr)^{\frac{n+4}{n-4}}+
{\frac{n+4}{n-4}}\sum_{j\ne
i}\bigl(\alpha_iP\delta_i \bigr)^{\frac{8}{n-4}}\alpha_j P\delta_j\\
 &+O\biggl(\sum_{j\ne i}P\delta_j ^{\frac{8}{n-4}}
P\delta_i\chi _{P\d_i\leq \sum_{j\neq i}P\d_j}
  +\sum_{j\ne i}P\delta_i ^{\frac{12-n}{n-4}}
 P\delta_j ^2 \chi_{P\d_j \leq P\d_i}
 +\sum_{k\ne j,k,j\ne i} P\delta_j ^{\frac{8}{n-4}}
 P\delta_k \biggr).\notag
\end{align}
Using Proposition \ref{13}, a computation similar to the one
performed in \cite{B1} and \cite{R} shows that\
\begin{align}\label{e:33}
 \biggl(P\d,\l \frac{\partial P\d}{\partial\l}\biggr)_2 & =\frac{n-4}{2}c_2
 \frac{H(a,a)}{\l^{n-4}}+O\biggl(\frac{1}{(\l d)^{n-2}}\biggr)\\
 \int K P\d ^{\frac{n+4}{n-4}}\l\frac{\partial P\d}{\partial\l} & =
-\frac{n-4}{n}c_3\frac{\D K(a)}{\l^2}+(n-4)c_2K(a)\frac{H(a,a)}{\l
^{n-4}}+ O\biggl(\frac{1}{\l^3}+\frac{1}{(\l
d)^{n-2}}\biggr).\notag
\end{align}
For $i\ne j$, we have
\begin{eqnarray}
\int_{\R^n}\d_i ^{\frac{n+4}{n-4}}\l_j\frac{\partial \d_j}
{\partial\l_j}=c_2\l_j\frac{\partial\e_{ij}}{\partial\l_j}+
O(\e_{ij}^{\frac{n-2}{n-4}}),
\end{eqnarray}
\begin{eqnarray}\label{e:35}
\biggl (P\d_j,\l_i\frac{\partial P\d_i}{\partial\l_i}\biggr)_2 =
c_2\biggl(\l_i\frac{\partial\e_{ij}}{\partial\l_i}+\frac{n-4}{2}
\frac{H(a_i,a_j)}{ (\l_i\l_j)^{(n-4)/2}}\biggr)+
O\biggl(\sum_{k=i,j} \frac{1}{(\l_kd_k)^{n-2}}+
\e_{ij}^{\frac{n-2}{n-4}}\biggr),
\end{eqnarray}
\begin{align}\label{e:36}
\int K& P\d_j^{\frac{n+4}{n-4}}\l_i\frac{\partial P\d_i}{
\partial\l_i} =K(a_j)\bigl(P\d_j,\l_i\frac{\partial P\d_i}{
\partial\l_i}\bigr)_2+O\biggl(\e_{ij}
(log\e_{ij}^{-1})^{\frac{n-4}{n}}\bigl(\frac{1}{\l_j}
+\frac{1}{(\l_jd_j)^4}\bigr)\biggr)\notag\\
  & +(\mbox{if }n\geq 8)
O\biggl(\e_{ij}^{\frac{n}{n-4}}log\e_{ij}^{-1}+
\frac{log(\l_jd_j)}{(\l_jd_j)^n}\biggr)+(\mbox{if }n<8) O\biggl
(\frac{\e_{ij}(log\e_{ij}^{-1})^{\frac{n-4}{n}}}
{(\l_jd_j)^{n-4}}\biggr),
\end{align}
\begin{align}\label{e:37}
 \int & K P\d_j \l_i\frac{\partial (P
\d_i)^{\frac{n+4}{n-4}}} {\partial\l_i}=
K(a_i)\bigl(P\d_j,\l_i\frac{\partial
P\d_i}{\partial\l_i}\bigr)_2+O\biggl(\e_{ij}
(log\e_{ij}^{-1})^{\frac{n-4}{n}}\bigl(\frac{1}{\l_i}
+\frac{1}{(\l_id_i)^4}\bigr)\biggr)\notag\\
 & +(\mbox{if }n\geq 8)
O\biggl(\e_{ij}^{\frac{n}{n-4}}log\e_{ij}^{-1}+
\frac{log(\l_id_i)}{(\l_id_i)^n}\biggr)+(\mbox{if }n<8)
O\biggl(\frac{\e_{ij}(log\e_{ij}^{-1})^{\frac{n-4}{n}}}
{(\l_id_i)^{n-4}}\biggr).
\end{align}
Now, it is easy to check
\begin{eqnarray}\label{e:38}
\mid\lambda_i{
\partial P\delta_i/\partial\lambda_i}\mid\leq c \delta_i, \quad
P\delta_k\leq \delta_k  \mbox{ and }
J(u)^{\frac{n}{n-4}}\a_j^{\frac{8}{n-4}}K(a_j)=1+o(1)\, \,
\forall\ j=1,...,p.
\end{eqnarray}
Combining \eqref{e:31},..., \eqref{e:38}, we easily derive our
proposition.
\end{pf}
\begin{pro}\label{p:34}
 For ${u=\sum_{i=1}^p \alpha_iP \delta_i}$ belonging to
$V(p,\varepsilon)$, we have the following expansion
\begin{align*}
\biggl(\partial J(u),\frac{1}{\l_i}\frac{\partial P\d_i}{\partial
a_i}\biggr)_2 & = 2J(u)\biggl[ -c_4\a_i ^{\frac{n+4}{n-4}}
J(u)^{\frac{n}{n-4}} \frac{\n K(a_i)}{\l_i}(1+o(1))\\
 & +
\frac{c_2}{2}\frac{\a_i}{ \l_i^{n-3}}\frac{\partial
H(a_i,a_i)}{\partial a_i}\bigl(1+o(1)\bigr)
 +O\biggl( \frac{1}{\l_i ^2}+\frac{1}{(\l_id_i)^{n-2}}
+\sum_{j\ne i}\e_{ij}\biggr)\biggr].
\end{align*}
We can improve this expansion and we obtain
\begin{align*}
\biggl(\partial & J(u), \frac{1}{\l_i}\frac{\partial P \d_i}{\partial
a_i}\biggr)_2= 2J(u)\biggl[ -c_4\a_i ^{\frac{n+4}{n-4}}
J(u)^{\frac{n}{n-4}} \frac{\n K(a_i)}{\l_i}(1+o(1))+
\frac{c_2}{2}\frac{\a_i}{ \l_i ^{n-3}}\frac{\partial
H(a_i,a_i)}{\partial a_i}\\
 & +c_2\sum_{j\ne i}\a_j\biggl(\frac{1}{\l_i}\frac{\partial
\e_{ij}}{\partial a_i} -\frac{1}{(\l_i \l_j)^{
\frac{n-4}{2}}}\frac{1}{\l_i}\frac{\partial H(a_i,a_j)}{\partial
a_i} \biggr)\biggl(1-J(u)^{\frac{n}{n-4}}\sum_{k=i,j}\a_k
^{\frac{8}{n-4}}K(a_k)\biggr)\biggr]\\
 & +O\biggl(\frac{1}{\l_i ^2}+\sum_{j\ne i} \l_j\mid
a_i-a_j\mid\e_{ij}^{\frac{n-1}{n-4}}\biggr)+o\biggl(\sum_k
\frac{1}{\l_k^2}+\frac{1}{(\l_kd_k)^{n-3}}+\sum_{k\ne j}
 \e_{kj}^{\frac{n-3}{n-4}}\biggr).
\end{align*}
\end{pro}
\begin{pf}
As in the proof of Proposition \ref{p:33}, we get \eqref{e:31} but
with ${\lambda_i{\partial P\delta_i/\partial\lambda_i}}$ changed
by ${{\lambda_i ^{-1}}{\partial P\delta_i/\partial a_i}}$.\\
Now, using Proposition \ref{13}, we observe (see \cite{B1} and
\cite{R})
\begin{align}\label{e:39}
\bigl(P\d,\frac{1}{\l}\frac{\partial P\d}{\partial a} \bigr)_2 &
=-\frac{c_2}{ 2\l^{n-3}}\frac{\partial H}{\partial a}(a,a)
+O\bigl(\frac{1}{(\l d)^{n-2}}\bigr), \\
\int K
P\d^{\frac{n+4}{n-4}}\frac{1}{\l}\frac{\partial P\d}{\partial a}
 & =-K(a)\frac{c_2}{\l^{n-3}}\frac{\partial H}{\partial a}(a,a)
+c_4\frac{\n K(a)}{\l}(1+o(1))+O\biggl(\frac{1}{\l^2}+\frac{1}{(\l
d)^{n-2}}\biggr)\notag
\end{align}
 where $c_4$ is a positive constant.\\
We also observe, for $i \ne j$
\begin{eqnarray}\label{o:32}
\int_{\R^n}\d_i ^{\frac{n+4}{n-4}}\frac{1}{\l_j} \frac{\partial
\d_j}{\partial
a_j}=c_2\frac{1}{\l_j}\frac{\partial\e_{ij}}{\partial a_j}+
O(\l_i|a_i-a_j|\e_{ij}^{\frac{n-1}{n-4}}),
\end{eqnarray}
\begin{align}\label{e:311}
\bigl(P\d_j,\frac{1}{ \l_i}\frac{\partial P\d_i}{\partial
a_i}\bigr)_2= & c_2\frac{1}{\l_i} \frac{\partial\e_{ij}}{\partial
a_i} -\frac{c_2}{ (\l_i\l_j)^{\frac{n-4}{2}}}\frac{1}{\l_i} \frac{\partial
H}{\partial a_i}(a_i,a_j)\\
 & +O\biggl(\sum_{k=i,j}\frac{1}{(\l_kd_k)^{n-2}}+
\e_{ij}^{\frac{n-1}{n-4}} \l_j\mid a_i-a_j\mid\biggr),\notag
\end{align}
\begin{align}\label{e:312}
\int K& P\d_j^{\frac{n+4}{n-4}} \frac{1}{\l_i}\frac{\partial
P\d_i}{\partial a_i}  =K(a_j)\bigl(P\d_j,\frac{1}{ \l_i}\frac{\partial
P\d_i}{\partial a_i}\bigr)_2+O\biggl(\e_{ij}
(log\e_{ij}^{-1})^{\frac{n-4}{n}}\bigl(\frac{1}{\l_j}
+\frac{1}{(\l_jd_j)^4}\bigr)\biggr)\notag\\
 & +(\mbox{if }n\geq 8)
O\biggl(\e_{ij}^{\frac{n}{n-4}}log\e_{ij}^{-1}+
\frac{log(\l_jd_j)}{(\l_jd_j)^n}\biggr)+(\mbox{if }n<8)
O\biggl(\frac{\e_{ij}(log\e_{ij}^{-1})^{\frac{n-4}{n}}}
{(\l_jd_j)^{n-4}}\biggr),
\end{align}
\begin{align}\label{e:313}
 \int K & P\d_j\frac {1}{\l_i}
\frac{\partial (P\d_i)^{\frac{n+4}{n-4}}}{\partial
a_i}=K(a_i)\bigl(P \d_j,\frac{1}{\l_i}\frac{\partial
P\delta_i}{\partial a_i}\bigr)_2 +O\biggl(\e_{ij}
(log\e_{ij}^{-1})^{\frac{n-4}{n}}\bigl(\frac{1}{\l_i}
+\frac{1}{(\l_id_i)^4}\bigr)\biggr)\notag\\
 & +(\mbox{if }n\geq 8)
O\biggl(\e_{ij}^{\frac{n}{n-4}}log\e_{ij}^{-1}+
\frac{log(\l_id_i)}{(\l_id_i)^n}\biggr)+(\mbox{if }n<8)
O\biggl(\frac{\e_{ij}(log\e_{ij}^{-1})^{\frac{n-4}{n}}}
{(\l_id_i)^{n-4}}\biggr).
\end{align}
Using \eqref{e:39},..., \eqref{e:313}, the proposition follows.
\end{pf}
\section{Proof of Theorem \ref{t:11} and its Corollary}
\mbox{}

First, we prove the following technical result which will be useful to prove the positivity of the solution that we will find.
\begin{pro}\label{p:40}
There exists a positive constant $\e_0$ such that, if $u \in H^2(\O )$ is a solution of the following equation
$$
\D^2 u = K |u| ^{\frac{8}{n-4}}u\mbox{ in }\O, \quad u=\D u=0
\mbox{ on }\partial\O
$$
and satisfying $|u^-|_{L^{\frac{2n}{n-4}}} < \e_0 $,
then $u$ has to be positive.
\end{pro}

\begin{pf}
First,  we observe that $K (u^-)^{\frac{n+4}{n-4}}\in L^{\frac{2n}{n+4}}$, where $u^-=\max (0,-u)$. \\
Now, let us introduce $w$ satisfying
\begin{eqnarray}\label{e:w}
\D^2 w=-K (u^-) ^{\frac{n+4}{n-4}} \mbox{ in }\O ,\quad w=\D w =0 \mbox{ on }\partial \O.
\end{eqnarray}
 Using a regularity argument, we derive that $w\in H^2\cap H_0^1(\O)$. Furthermore, the maximum principle  implies that $w\leq 0$. Now, multiplying  equation \eqref{e:w} by $w$ and integrating on $\O$, we derive that
\begin{eqnarray}
||w||_2^2=\int_\O \D^2 w .w=-\int_\O K(u^-)^{\frac{n+4}{n-4}}w\leq c||w||_2|u^-|_{L^{\frac{2n}{n-4}}}^{\frac{n+4}{n-4}}.
\end{eqnarray}
Thus, either $||w||_2=0$ and it follows that $u^-=0$ or $||w||_2\ne 0$ and therefore
\begin{eqnarray}
||w||_2\leq c|u^-|_{L^{\frac{2n}{n-4}}}^{\frac{n+4}{n-4}}.
\end{eqnarray}
On the other hand, we have
\begin{eqnarray}
\int_\O \D^2 w.u =\int_\O K(u^-)^{\frac{2n}{n-4}}\geq c|u^-|_{L^{\frac{2n}{n-4}}}^{\frac{2n}{n-4}}.
\end{eqnarray}
Furthermore we obtain
\begin{align}
\int_\O \D^2w.u & =\int_\O w.\D^2u=\int_\O K |u|^{\frac{8}{n-4}}u w=-\int_{u\leq 0}K(u^-)^{\frac{n+4}{n-4}}w +\int_{u\geq 0}K(u^+)^{\frac{n+4}{n-4}}w \\
 & \leq \int_{u\leq 0}-K(u^-)^{\frac{n+4}{n-4}}w=\int_\O-K (u^-)^{\frac{n+4}{n-4}}w=\int_\O \D^2 w.w
 =||w||_2^2.
\end{align}
Thus,
\begin{eqnarray}
|u^-|_{L^{\frac{2n}{n-4}}}^{\frac{2n}{n-4}}\leq c||w||_2^2\leq c |u^-|_{L^{\frac{2n}{n-4}}}^{\frac{2(n+4)}{n-4}}.
\end{eqnarray}
Thus, for $|u^-|_{L^{\frac{2n}{n-4}}}$ small enough, we derive a contradiction and therefore the case $||w||_2 \ne 0$ cannot occur, so $||w||_2$ has to be equal to zero and therefore $u^-=0$. Thus the  result follows.
\end{pf}

Now, we provide the characterization of the critical points at infinity of
$J$ in the case where we have only  one mass. We recall that
the critical points at infinity are the orbits of the gradient
flow of $J$ which remain in $V(p,\e (s))$, where $\e (s)$ is some
 function such that $\e(s)$ tends to zero when $s$ tends to $+\infty$,
see \cite{B1}.
\begin{pro}\label{p:41}
Let $n\geq 7$ and assume that $(A_0)$ holds. Then there exists a pseudogradient $Y_1$ such that the
following holds:\\
 there exists a constant $c>0$ independent of $u=\a\d_{(a,\l)}\in
 V(1,\e)$ such that
  $$(-\partial J(u),Y_1)_2\geq c\biggl(\frac{1}{\l^2}+\frac{\mid\n
  K(a)\mid}{\l}+\frac{1}{(\l d)^{n-3}}\biggr)\leqno{1)}$$
 $$(-\partial J(u+\ov{v}),Y_1+\frac{\partial\ov{v}}{\partial(\a,a,\l)}
 (Y_1))_2\geq c\biggl(\frac{1}{\l^2}+\frac{\mid\n
  K(a)\mid}{\l}+\frac{1}{(\l d)^{n-3}}\biggr)\leqno{2)}$$
 $3)$ $ Y_1$ is bounded and the only case where $\l$ increases along
 $Y_1$ is when $a$ is close to a critical point $y$ of $K$ with $-\D
 K(y)>0$. Furthermore the distance to the boundary only increases
 if it is small enough.
 \end{pro}
 \begin{pf}
Using $(A_0)$ and the fact that the boundary of $\O$ is a compact
set, then there exist two positive constants $c$ and $d_0$ such
that for each $x$ satisfying $d_x\leq d_0$ we have $\n
K(x).\nu_x<-c$ where $\nu_x$ is the outward normal to
$\O_{d_x}=\{z\in \O/d_z= d(z,\partial\O) > d_x\}$.
The construction will depend on $a$ and $\l$. We distinguish three cases:\\
 {\it 1st case: } If $a$ is near the boundary, that is $d_a\leq d_0$,
 we define
 $$W_1=-\frac{1}{\l}\frac{\partial P\d_{(a,\l)}}{\partial a}\nu_a.$$
 {\it 2nd case:} If $d_a\geq d_0$ and $|\n K(a)|\geq C_2/\l$
 where $C_2$ is a large positive constant. In this case, we define
 $$W_2=\frac{1}{\l}\frac{\partial P\d_{(a,\l)}}{\partial a}
 \frac{\n K(a)}{|\n K(a)|}.$$
 {\it 3rd case: } If $|\n K(a)|\leq 2C_2/\l$, thus $a$ is near
a critical point $y$ of $K$. Then we define
 $$W_3=(sign(-\D K(y)))\l\frac{\partial P\d_{(a,\l)}}{\partial\l}.
$$
 In all cases, using Propositions \ref{p:33} and
 \ref{p:34}, we derive that
 $$\bigl(-\partial J(u),W_i\bigr)_2\geq c\biggl(\frac{1}{\l^2}
 +\frac{1}{(\l d)^{n-3}}+\frac{\mid\n K(a)\mid}{\l}\biggr).$$
 The pseudogradient $Y_1$ will be a convex combination of $W_1$,
$W_2$ and $W_3$. Thus the proof of claim 1) is completed. The proof of claim
2) follows from the estimate of $\ov{v}$ as in \cite{B2} and
\cite{BCCH}. The proof of claim 3) follows from the construction of the
vector field $Y_1$.
 \end{pf}
\begin{pro}\label{p:42}
Assume that $J$ does not have any critical points in $\Sig^+$ and
assume that $(A_0)$ and $(A_2)$ hold. Then the only critical
points at infinity of $J$ in $V(1,\e )$, for $\e$ small enough,
correspond to
 $P\d_{(y,+\infty)}$ where $y$ is a critical point of $K$ with
 $-\D K(y)>0$ if $n\geq 7$ and with $-\D K(y)/(60 K(y)) + H(y,y)
 > 0$ if $n=6$.  Moreover, such a critical point at infinity has
 a Morse index equal
 to $n-index(K,y)$.
\end{pro}
\begin{pf}
First, we recall that the $6$-dimension case of such a Proposition
has already been proved in \cite{BH}, so we need to prove our
result
for $n\geq 7$. \\
Now, from Proposition \ref{p:41}, we know that the only region
where $\l$ increases along the pseudogradient $Y_1$, defined in
Proposition \ref{p:41}, is the region where $a$ is near a critical
point $y$ of $K$ with $-\D K(y) > 0$. Arguing as in \cite{B2} and
\cite{BCCH}, we can easily derive
from Proposition \ref{p:41}, the following normal form :\\
if $a$ is  near a critical point $y$ of $K$ with $-\D K(y)>0$,
we can find a change of variables
$(a,\l)\longrightarrow
(\bar{a},\bar{\l})$ such that
\begin{eqnarray}\label{mok}
J(P{\d}_{(a,\l)} + \bar{v})=\Psi
(\bar{a},\bar{\l}):=
\frac{S_n^{4/n}}{K(\bar{a})^{(n-4)/n}}\left(1-\frac{(c-\eta)}
{\bar\l^2}\frac{\D K(y)}{K(y)^{n/4}}\right),
\end{eqnarray}
where $c$ is a constant which depends only on $n$ and $\eta$
 is a small positive constant.\\
This  yields a split of variables $a$ and $\l$, thus it
follows that if $a=y$, only $\l$ can move. In order to
decrease the functional $J$, we have to increase $\l$,
thus we find a critical point at infinity only in this
case and our result follows.
\end{pf}\\
Now, we are ready to prove Theorem \ref{t:11} and its corollary.\\
\begin{pfn}{\bf Theorem \ref{t:11}}
 Arguing by
contradiction, we suppose that $J$ has no critical points in
$\Sig^+$. It follows from Proposition \ref{p:31} and Proposition \ref{p:42}, that under the
assumptions of Theorem \ref{t:11}, the critical points at infinity
of $J$ under the level $c_1=
(S_n)^{\frac{4}{n}}(K(y_l))^{\frac{4-n}{n}} + \e $ , for $\e$
small enough, are in one to one correspondence with the  critical
points of $K$ $y_0$, $y_1$, ..., $y_l$. The unstable manifold at
infinity of such  critical points at infinity, $W_u(y_0)_\infty$,
..., $W_u(y_l)_\infty$ can be described, using \eqref{mok}, as the
product of $W_s(y_0)$, ..., $W_s(y_l)$ (for a pseudogradient of
$K$ ) by $[A, +\infty [$
 domain of the variable $\l$, for some positive
number $A$ large enough.\\
Let $\eta$ be a small positive constant and let \be\label{**}
V_\eta (\Sig ^+)=\{u\in \Sig / J(u)^{\frac{2n-4}{n-4}}e^{2J(u)}
|u^-|^{\frac{8}{n-4}}_{L^{\frac{2n}{n-4}}} < \eta \}. \ee Since
$J$ has no critical points in $\Sig^+$, it follows that $
J_{c_1}=\{u\in V_{\eta}( \Sig ^+) / J(u) \leq c_1 \}$ retracts by
deformation on $X_\infty = \cup _{0\leq j\leq l}W_u(y_j)_\infty$
(see Sections 7 and 8 of \cite{BR}) which can be parametrized as
we said before by $X
\times [A, +\infty[$.\\
On the other hand, we have $X_\infty$ is contractible in
$J_{c_2+\e}$, where
$c_2=(S_n)^{\frac{4}{n}}\bar{c}^{\frac{4-n}{n}}$. Indeed, from
$(A_4)$, it follows that there exists a contraction $ h :[0,1]
\times X \to K^{\bar{c}}$, $h$ continuous, such that for any $a\in
X$, $h(0,a)=a$ and $h(1,a)=a_0\in X$. Such a contraction gives
rise to the following contraction $\tilde{h} : [0,1]\times
X_\infty \to V_{\eta}(\Sig ^+)$ defined by
$$[0,1] \times X \times \left[A,\right.+\infty\left[ \right.\ni
(t,a,\l  ) \longmapsto P\d _{(h(t,a),\l )} + \bar{v} \in
V_{\eta}(\Sigma ^+).
$$
In fact, $\tilde{h}$ is continuous and it satisfies
$\tilde{h}(0,a,\l) = P\d _{(a, \l) } +\bar{v} \in X_\infty$ and
$\tilde{h}(1,a,\l )= P\d _{(a_0,\l)} +\bar{v}$. \\
Now, using Proposition \ref{p:31}, we deduce that
$$
J(P\d _{(h(t,a), \l) } + \bar{v}) =
({S_n})^{\frac{4}{n}}(K(h(t,a)))^{\frac{4-n}{n}}\left(1+O(A
^{-2})\right),
$$
where $K(h(t,a)) \geq \bar{c} $ by construction.\\
Therefore such a contraction is performed under $c_2 +\e$, for $A$
large enough, so $X_\infty$ is contractible in $J_{c_2+\e }$.\\
In addition, choosing $c_0$ small enough, we see that there is no
critical point at infinity for $J$ between the levels $c_2+\e$ and
$c_1$, thus  $J_{c_2+\e }$ retracts by deformation on $J_{c_1}$,
which retracts by deformation on $X_\infty$, therefore $X_\infty$
is contractible leading to the contractibility of $X$, which is in
contradiction with assumption $(A_3)$. Hence $J$ has a critical
point in $V_\eta (\Sig ^+ )$. Using Proposition \ref{p:40}, we
derive that such a critical point is positive. Therefore our
theorem follows.
 \end{pfn}

Now, we give the proof of Corollary \ref{c:12}.

\begin{pfn}{\bf Corollary \ref{c:12}}
Arguing by contradiction, we may assume that the Morse index of
the
solution provided by Theorem \ref{t:11} is $\leq m-1$.\\
Perturbing, if necessary $J$, we may assume that all the critical
points of $J$ are nondegenerate and have their Morse index $\leq
m-1$. Such critical points do not change the homological group in
dimension $m$ of level sets of $J$.\\
Since $X_\infty$ defines a homological class in dimension $m$
which is nontrivial in $J_{c_1}$, but trivial in $J_{c_2+\e}$, our
result follows.
\end{pfn}
\section{Proof of Theorem \ref{t:13}}
\mbox{}
 Arguing by
contradiction, we suppose that $J$ has no critical points in $V_\eta(\Sig ^+)$ defined by \eqref{**}. We
denote by $z_1,...,z_r$ the critical points of $K$ among of $ y_i$
$(1\leq i\leq l )$, where
$$
 -\D K(z_j) \leq 0  \quad
(1\leq j\leq r).
$$
The idea of the Proof of Theorem \ref{t:13} is to perturb the
function $K$ in the $C^1$ sense in some  neighborhoods of
$z_1,...,z_r$ such that the new function $\tilde{K}$ has the same
critical points with the same  Morse indices but satisfying $-\D
\tilde{K}(z_j) > 0$ for $1\leq j\leq r$. Notice that the new
$\tilde{X}$ corresponding to $\tilde{K}$, defined in assumption
$(A_3)$, is also not contractible and its homology group in
dimension $m$ is nontrivial.

Under the level $2^{4/n}S_n^{4/n}(K(y_0))^{(4-n)/n}$, the
associated functional $\tilde{J}$ is close to the functional $J$
in the $C^1$ sense. Under the level $c_2 +\e$, where $c_2$ is
defined in the proof of Theorem \ref{t:11}, the functional
$\tilde{J}$ may have other critical points, however a careful
choice of $\tilde{K}$ ensures that all these critical points have
Morse indices less than $m - 2$ (see Proposition \ref{p:51}
below), and so they do not change the homology in dimension $m$,
therefore the arguments used in the Proof of Theorem \ref{t:11}
lead to a contradiction. It follows that Theorem \ref{t:13} will
be a corollary of the following Proposition:
\begin{pro}\label{p:51}
Assume that $J$ has no critical points in $V_\eta(\Sig ^+)$. We
can choose $\tilde{K}$ close to $K$ in the $C^1$ sense such that
$\tilde{K}$ has the same critical points with the same  Morse
indices and such that:
\begin{align*}
i)\quad & -\D \tilde{K}(z_j) > 0 \quad\mbox{ for}\quad  1\leq j\leq r,\\
ii)\quad &  -\D \tilde{K}(y) > 0 \quad\mbox{ for}\quad
y\in\{y_0,..., y_l\} \diagdown \{z_1,..., z_r\},\\
iii)\quad & -\D \tilde{K}(y_i) < 0 \quad\mbox{ for}\quad  l+1\leq i\leq s,\\
iv)\quad &\mbox{ if } \tilde{J} \mbox{ has critical points under the level } c_2 +\e, \mbox{ then their Morse}\\
\quad &\mbox{  indices are less than } m-2, \mbox{ where } m  \mbox{ is defined in assumption } (A_3),\\
v)\quad & \mbox{ the new } \tilde{X} \mbox{ corresponding to } \tilde{K}, \mbox{ defined in assumption }(A_3), \mbox{ is also}\\
\quad & \mbox{  not contractible and its homology group in
dimension } m \mbox{ is nontrivial}.
\end{align*}
\end{pro}

Next, we are going to prove Proposition \ref{p:51}. For this purpose,
we need the following lemmas.
\begin{lem}\label{l:52}
 Let $z_0$ be a point of $\O$ such that $ d(z_0,\partial\O ) \geq c_0
 >0$ and let $\pi$ be the orthogonal projection (with respect
 to the scalar inner $(u,v)_2= \int_\O \D u\D v$) onto\\  $E^\bot =
 \mbox{ Vect }\left(P\d_{(z_0,\l )}, \l^{-1}\partial
 P\d_{(z_0,\l )}/\partial z, \l \partial
 P\d_{(z_0,\l )}/\partial\l\right)$. Then, we have the
 following estimates
$$
(i)\quad ||J'(P\d_{(z_0,\l )})|| = O\left(\frac{1}{\l}\right);
\quad (ii)\quad || \frac{\partial\pi}{\partial z}||= O(\l);
\quad (iii)\quad || \frac{\partial ^2 \pi}{\partial ^2 z}||=
 O(\l ^2).
$$
\end{lem}
\begin{pf}
The proof of claim $(i)$ is easy, so we will omit it. Now, we
prove claim $(ii)$. Let $\var \in \{P\d_{(z_0,\l )}, \l^{-1}
\partial P\d_{(z_0,\l )}/\partial z, \l \partial P\d_{(z_0,\l
)}/\partial\l\}$. We then have $\pi\var =\var$,
 therefore
$$
\frac{\partial \pi}{\partial z}(\var )= \frac{\partial \var}
{\partial z} - \pi\frac{\partial \var}{\partial z},
$$
thus $||\frac{\partial \pi}{\partial z}(\var )||= O(\l )$.\\
Now, for $v\in E$, we have $\pi v=0$, thus
$$
\frac{\partial \pi}{\partial z}v = - \pi\frac{\partial v}
{\partial z}= \sum_{i=1}^3a_i\var _i ,
$$
where $\var _1 =P\d_{(z_0,\l )}$, $\var _2=  \l^{-1}\partial
P\d_{(z_0,\l )}/\partial z$, $\var _3= \l \partial P\d_{
(z_0,\l )}/\partial\l$.\\
But, we have
$$
a_i||\var _i||^2 = (\frac{\partial v}{\partial z}, \var _i)_2= -
(v, \frac{\partial\var_i}{\partial z})_2 = O(\l ||v||).
$$
Thus claim $(ii)$ follows.\\
 In the same way, claim $(iii)$ follows and hence the proof of
 our lemma is completed.
\end{pf}
\begin{lem}\label{l:53}
Let $z_0$ be a point of $\O$ close to a critical point of $K$
such that $d(z_0,\partial\O )\geq c_0>0$. Let $\bar{v}=
\bar{v}(z_0,\a , \l )\in E$ defined in Proposition \ref{p:32}.
Then, we have the following estimates
$$
(i) \qquad ||\bar{v}||= o(\frac{1}{\l}), \qquad \quad (ii) \qquad
||\frac{\partial \bar{v}}{\partial z}|| = o(1).
$$
\end{lem}
\begin{pf}
We notice that Claim $(i)$ follows from Proposition \ref{p:32}.
Then, we  need only to show that Claim $(ii)$ is true. We know that
$\bar{v}$ satisfies
$$
A\bar{v}= f + O\left(||\bar{v}||^{(n+4)/(n-4)}\right)\quad \mbox{and }
 \quad
\frac{\partial A}{\partial z} \bar{v} +A \frac{\partial
\bar{v}}{\partial z} = \frac{\partial f}{\partial z} +
O\biggl(||\bar{v}||^{8/(n-4)} |\frac{\partial \bar{v}}{\partial
z}|\biggr),
$$
where $A$ is the operator associated to the quadratic form $Q$
defined on $E$ ($Q$ and $f$ are defined in Proposition \ref{p:31}).\\
Then, we have
$$
A\biggl(\frac{\partial\bar{v}}{\partial z} -\pi(\frac{\partial
\bar{v}}{\partial z})\biggr) = \frac{\partial f}{\partial z} -
\frac{\partial A}{\partial z} \bar{v}- A\pi(\frac{\partial
\bar{v}}{\partial z}) +  O\biggl(||\bar{v}||^{8/(n-4)}
|\frac{\partial \bar{v}}{\partial z}|\biggr).
$$
Since $Q$ is a positive quadratic form on $E$ (see \cite{BE1}), we
then derive
$$
||\frac{\partial\bar{v}}{\partial z} -\pi(\frac{\partial
\bar{v}}{\partial z})|| \leq C\biggl(||\frac{\partial f}{
\partial z}|| + ||\frac{\partial A}{\partial z}||||
\bar{v}|| +||\pi(\frac{\partial \bar{v}}{\partial z})||+ ||
\bar{v}||^{\frac{8}{n-4}}||\frac{\partial \bar{v}}{
\partial z}||\biggr).
$$
Now, we estimate each term of the right hand-side in
the above estimate. First, it is easy to see $||\frac{
\partial A}{\partial z}||=O(\l )$. Therefore, using (i),
we obtain $||\frac{\partial A}{\partial z}|| || \bar{v}||
=o(1)$. Secondly, we have
\begin{align}
(\frac{\partial f}{\partial z},v)_2 & =c\int KP
\d_{(z_0,\l)}^{\frac{8}{n-4}}\frac{\partial P\d}{\partial z}v=c\n
K(z_0)\int d(z_0,x)\d ^{\frac{8}{n-4}}\frac{\partial \d }{\partial
z}v \notag\\
 & +O\biggl(\int d^2(x,z_0)\d^{\frac{n+4}{n-4}}
\l |v|\biggr)+ O\left(\int_\O \d^{8/(n-4)}\var |v| +\int_\O
\d^{8/(n-4)} |\frac{\partial\var}{\partial z}| |v|\right)\notag \\
 & \leq c ||v||(|\n K(z_0)|+\frac{1}{\l}),
\end{align}
where $\var = \d -P\d$.\\
Since $z_0$ is close to a critical point of $K$, we derive that
$||\frac{\partial f}{\partial z}||=o(1)$.\\
For the  term $||\pi (\frac{\partial \ov{v}}{\partial z})||$, we
have, since $\ov{v}\in E$
\begin{align*}
\left(\frac{\partial \ov{v}}{\partial z},\d_{(z_0,\l)}\right)_2= &
-
\left(\ov{v},\frac{\partial \d_{(z_0,\l)}}{\partial z}\right)_2=0\\
\left(\frac{\partial \ov{v}}{\partial z},\l\frac{ \partial
\d_{(z_0,\l)}}{{\partial \l}}\right)_2= & -  \left(\ov{v},\l
\frac{\partial^2 \d_{(z_0,\l)}}{\partial\l\partial z}\right)_2
=O(\l||\ov{v}||)=o(1)
\end{align*}
In the same way, we have
$$
\left(\frac{\partial \ov{v}}{\partial
z},\frac{1}{\l}\frac{\partial P\d}{{\partial z}}\right)_2=o(1)
$$
Therefore $||\pi(\frac{\partial \ov{v}}{\partial z})||=o(1)$.
Now, using the following inequality
$$
||\frac{\partial \ov{v}}{\partial z}||\leq ||\frac{\partial
\ov{v}}{\partial z} -\pi (\frac{\partial \ov{v}}{\partial z})||
 + ||\pi (\frac{\partial \ov{v}}{\partial z})||,
$$
we easily derive our claim and our lemma follows.
\end{pf}\\
We are now able to prove Proposition \ref{p:51}.\\
\begin{pfn}{Proposition \ref{p:51}}
We suppose that $J$ has no critical points in $V_\eta (\Sig^+)$
and we perturb  the function $K$ only in some neighborhoods of
$z_1,...,z_r$, therefore Claims $ii)$ and $iii)$ follow from
assumption $(A_2')$. We observe that under the level $c_2+\e$ and
outside $V(1,\e_0)$, we have $|\partial J|> c > 0$. If $\wtilde K$
is close to $K$ in the  $C^1$-sense, then  $\wtilde J$ is close to
$J$ in the  $C^1$-sense, and therefore $|\partial \wtilde J|> c/2$
in this region. Thus, a critical point $u_0$ of $\wtilde J$ under
the level $c_2+\e$ has to be in $V(1,\e_0)$. Therefore, we can
write $u_0=P\d_{(z_0,\l)}+\ov{v}$.\\
 Next we will prove the following Claim \\
{\bf Claim:} $z_0$ has to be near a critical point $z_i$ of $K$,
 $1\leq i\leq r$ (recall that  $z_i$'s satisfy $\D K(z_i) \geq 0$).\\
To prove our Claim, we will prove in the first step that $d_{z_0}
 := d(z_0,\partial\O ) \geq c_0>0$. For this fact, arguing by
 contradiction, we assume that $d_{z_0} \to 0$. Thus, we have
\begin{eqnarray}
\frac{\partial K}{\partial\nu}(z_0) < -c<0 \qquad \mbox{and}
\qquad \frac{\partial H}{\partial\nu}(z_0,z_0) \thicksim
\frac{c}{d_{z_0}^{n-3}}
\end{eqnarray}
(the proof of the last fact is similar to the corresponding statement for the Laplacian operator in \cite{R}).\\
Using Propositions \ref{p:32} and \ref{p:34}, we obtain
$$
0=\left(\partial \wtilde J(u_0),\frac{1}{\l}\frac{\partial
P\d}{\partial z}\right)_2. \nu > \frac{c}{\l} + \frac{c}{(\l\
d_{z_0})^{n-3}} > 0
$$
Thus, we derive a contradiction and therefore $z_0$ has to satisfy
$d_{z_0} \geq c_0>0$.\\
 Now, also using Propositions \ref{p:32} and  \ref{p:34}, we derive that
$$
0=\left(\partial \wtilde J (u_0),\frac{1}{\l}\frac{\partial
P\d}{\partial z}\right)_2= c\frac{\n \wtilde
K(z_0)}{\l}+o(\frac{1}{\l}),
$$
thus $z_0$ has to be close to $y_i$ where $i\in \{0,...,s\}$.\\
 We also have, by Propositions \ref{p:32} and  \ref{p:34}
\begin{eqnarray}\label{e:51}
0=\left(\partial \wtilde J (u_0),\l\frac{\partial P\d}{\partial
\l}\right)_2=c\frac{\D \wtilde K(z_0)}{\l^2}+o(\frac{1}{\l^2})
\end{eqnarray}
In the neighborhood of $y_i$ with $i\in\{k/-\D K(y_k) >0\}\cup
\{l+1,...,s\}$, $\wtilde K\equiv K$ and therefore
$|\D\wtilde K| >c >0$ in this neighborhood. Thus \eqref{e:51}
implies that $z_0$ has to be near $z_i$ with $1\leq i\leq r$.
Thus our Claim is proved.\\
In the sequel, we assume that $\d=\d_{(z_0,\l)}$ satisfies
$||\d||=1$, and thus
$\D^2\d=S_n^{\frac{4}{n-4}}\d^{\frac{n+4}{n-4}}$. We also assume
that $|D^2 \wtilde K| \leq c(1+|D^2 K|)$,
where $c$ is a fixed positive constant. \\
Let $u_0=P \d_{(z_0,\l)}+\ov{v}$ be a critical point of $\wtilde
J$. In order to compute the Morse index of $\wtilde J$ at $u_0$,
we need to compute  $\frac{\partial^2}{\partial z^2}\wtilde
J(P\d_{(z,\l)}+\ov{v})_{|z=z_0}$.\\
 We observe that
 $$\frac{\partial}{\partial z}\wtilde J(P\d_{(z,\l)}+\ov{v}) =
\wtilde J'(P \d_{(z,\l)}+\ov{v})\frac{\partial}{\partial z}
(P\d_{(z,\l)}+\ov{v}) = \wtilde J'(P\d_{(z,\l)}+\ov{v}) \pi
 \frac{\partial}{\partial z}(P\d_{(z,\l)}+\ov{v}))
$$
and
\begin{align}\label{e:52}
\frac{\partial^2}{\partial z^2}\wtilde J(P\d_{(z,\l)}+\ov{v}) = &
\wtilde J''(P\d_{(z,\l)}+\ov{v})\frac{\partial}{\partial z}(P
\d_{(z,\l)}+\ov{v}) \pi
 (\frac{\partial}{\partial z}(P\d_{(z,\l)}+\ov{v}))\\
 & +
 \wtilde J'(P\d_{(z,\l)}+\ov{v})\frac{\partial}{\partial z}\biggl(
\pi (\frac{\partial}{\partial
z}(P\d_{(z,\l)}+\ov{v}))\biggr).\notag
\end{align}
For $z=z_0$, we have $\wtilde J'(P\d_{(z,\l)}+\ov{v})=0$. We will
estimate each term of the right hand-side of \eqref{e:52}. First,
we have by Lemma \ref{l:53}
$$
\wtilde J''(P\d_{(z,\l)}+\ov{v})\frac{\partial \ov{v}}{\partial z}
\pi (\frac{\partial \ov{v}}{\partial z})=o(1).
$$
Secondly, we compute
$$
T=\wtilde J''(P\d_{(z,\l)}+\ov{v})\frac{\partial P \d}{\partial z}
\pi
 \frac{\partial \ov{v}}{\partial z}=c\bigg[\left(\frac{\partial P\d}{\partial z},
\pi \frac{\partial \ov{v}}{\partial
z}\right)-\frac{n+4}{n-4}\wtilde J(u_0)^{\frac{n}{n-4}}\int
\tilde{K}(P \d+\ov{v})^{\frac{8}{n-4}}\frac{\partial P\d}{\partial
z} \pi \frac{\partial \ov{v}}{\partial z}\bigg]
$$
According to Proposition \ref{p:31}, we have
\begin{eqnarray}\label{mokh}
\wtilde J(P\d+\ov{v})=\frac{S_n^{4/n}}{\wtilde
K(z)^{\frac{n-4}{n}}}+O(\frac{||\ov{v}||}{\l}+\frac{1}{\l ^2}).
\end{eqnarray}
Thus
\begin{align*}
T & = c\bigg[\left(\frac{\partial P\d}{\partial z}, \pi
\frac{\partial \ov{v}}{\partial
z}\right)_2-\frac{n+4}{n-4}S_n^{\frac{4}{n-4}}\int
\frac{\tilde{K}}{\tilde{K}(z)}P\d ^{\frac{8}{n-4}}\frac{\partial
P\d}{\partial z} \pi
 \frac{\partial \ov{v}}{\partial z}\bigg]\\
 & +O\left(\int \left(
\d^{\frac{12-n}{n-4}}|\ov{v}|+|\ov{v}|^{\frac{8}{n-4}}\chi_{P
\d\leq |\ov{v}|}\right)|\frac{\partial P\d}{\partial z}|| \pi
 \frac{\partial \ov{v}}{\partial z}|\right)+o(1)\\
& = c\frac{n+4}{n-4}S_n^{\frac{4}{n-4}}\int \left(1-
\frac{\tilde{K}}{\tilde{K}(z)}\right)\d
^{\frac{8}{n-4}}\frac{\partial \d}{\partial z} \pi (\frac{\partial
\ov{v}}{\partial z})+ O\biggl(\l||\ov{v}||||\frac{\partial
\ov{v}}{\partial
z}||+\l||\ov{v}||^{\frac{n+4}{n-4}}||\frac{\partial
\ov{v}}{\partial z}||\biggr)+o(1)\\
 & =o(1).
\end{align*}
Thus \eqref{e:52} becomes
\begin{align*}
\frac{\partial^2}{\partial z^2}\wtilde J(P\d_{(z,\l)} & +\ov{v}) =
\wtilde J''(P\d_{(z,\l)}+\ov{v})\frac{\partial P\d}{\partial
z}\left( \frac{\partial P\d}{\partial z}
+\frac{\partial\ov{v}}{\partial z}\right)+\wtilde
J'(P\d_{(z,\l)}+\ov{v})
\frac{\partial ^2 P\d}{\partial z^2}+o(1)\\
 & =2\wtilde J(u_0)\bigg[\left(\frac{\partial P\d}{\partial z}+
\frac{\partial \ov{v}}{\partial z},\frac{\partial P\d}{\partial
z}\right)_2+\left(P\d +\ov{v},\frac{\partial^2 P\d}{\partial z^2}\right)_2\\
 & -\wtilde
J(u_0)^{\frac{n}{n-4}}\frac{n+4}{n-4}\biggr(\int
K(P\d+\ov{v})^{\frac{8}{n-4}}(\frac{\partial P\d}{\partial z})^2
+\int K(P\d+\ov{v})^{\frac{8}{n-4}}\frac{\partial
P\d}{\partial z}\frac{\partial \ov{v}}{\partial z}\biggr)\\
 & -\wtilde
J(u_0)^{\frac{n}{n-4}} \int K(P\d+\ov{v})^{\frac{n+4}{n-4}}
\frac{\partial^2 P\d}{\partial z^2}\bigg]+o(1)\\
&= 2\wtilde J(u_0)\bigg[\left(\frac{\partial P\d}{\partial z}+
\frac{\partial \ov{v}}{\partial z},\frac{\partial P\d}{\partial
z}\right)_2+\left(P\d +\ov{v},\frac{\partial^2 P\d}{\partial z^2}\right)_2\\
 & -\frac{n+4}{n-4}\wtilde
J(u_0)^{\frac{n}{n-4}}\biggr(\int K
P\d^{\frac{8}{n-4}}(\frac{\partial P\d}{\partial z})^2
+\frac{8}{n-4}\int K P\d ^{\frac{12-n}{n-4}}\ov{v}(\frac{\partial
P\d}{\partial z})^2\\
 &  +\int KP\d ^{\frac{8}{n-4}}\frac{\partial P\d}{\partial
z}\frac{\partial \ov{v}}{\partial z}
 + \frac{n-4}{n+4}\int
K P\d ^{\frac{n+4}{n-4}}\frac{\partial^2 P\d}{\partial z^2}+\int K
P\d ^{\frac{8}{n-4}}\ov{v}\frac{\partial^2 P\d}{\partial z^2}
\biggr)\bigg]+o(1).
\end{align*}
Using \eqref{mokh} and  Proposition \ref{13}, we derive that
\begin{align*}
\frac{\partial^2}{\partial z^2} & \wtilde J(P\d_{(z,\l)} +\ov{v})
=2 \wtilde J(u_0)\bigg[S_n^{\frac{4}{n-4}}
\biggl(\int\frac{n+4}{n-4} \d^{\frac{8}{n-4}}(\frac{\partial
\d}{\partial z})^2+\d^{\frac{n+4}{n-4}}\frac{\partial^2
\d}{\partial z^2}\biggr)\\
 & -\wtilde
J(u_0)^{\frac{n}{n-4}}\biggr(\frac{n+4}{n-4}\int K
\d^{\frac{8}{n-4}}(\frac{\partial \d}{\partial z})^2 +\int K \d
^{\frac{n+4}{n-4}}\frac{\partial^2 \d}{\partial z^2}\biggr)\\
 & + S_n^{\frac{4}{n-4}}\frac{n+4}{n-4}\biggl(\int
\bigl(1-\frac{K}{K(z)}\bigr)\d ^{\frac{8}{n-4}}\frac{\partial
\d}{\partial z}\frac{\partial \ov{v}}{\partial z}+\int
\bigl(1-\frac{K}{K(z)}\bigr)\d ^{\frac{8}{n-4}}\frac{\partial ^2
\d}{\partial z^2} \ov{v}\\
 & +\frac{8}{n-4}\int
\bigl(1-\frac{K}{K(z)} \bigr)\d ^{\frac{12-n}{n-4}}
(\frac{\partial \d}{\partial z})^2 \ov{v}\biggr)\bigg]+o(1)\\
 & =2\wtilde
J(u_0)\bigg[S_n^{\frac{4}{n-4}}\frac{\partial}{\partial
z}\biggl(\int_{R^n}\d^{\frac{n+4}{n-4}}\frac{\partial \d}{\partial
z}\biggr)-\wtilde J(u_0)^{\frac{n}{n-4}}\int_{\O}K
\frac{\partial^2\d^{\frac{2n}{n-4}}}{\partial z^2}\bigg]+o(1).
\end{align*}
Thus
$$\frac{\partial^2}{\partial z^2}\wtilde J(P\d_{(z,\l)}
+\ov{v})_{|_{z=z_0}} =-c D^2 K(z_0)+o(1),
$$
where $c$ is a positive constant.\\
 Therefore, taking account of the $\l$-space, we derive that
$$\mbox{index}(\tilde{J},u_0)\leq n- \mbox{index}(K,z_0)+1\leq m-2.$$
Then  Claims $i)$ and $iv)$ of Proposition \ref{p:51} follow.\\
On the other hand, according to  assumption $(A_2')$ we have
$$
n-m+3 \leq \mbox{index} (K,z_j)= \mbox{index} (\tilde{K}, z_j)
\quad\mbox{for}\quad 1\leq j\leq r.
$$
Thus, for any pseudogradient of $\tilde{K}$, the dimension of the
stable manifold of $z_j$ is less than $m-3$. Note that our
perturbation changes the pseudogradient $Z$ to $\tilde{Z}$, but
only in some neighborhoods of $z_1,..., z_r$. Therefore the stable
manifolds of $y_i$, for $i \not\in \{1,...,r\}$, remain unchanged.
Since the dimension of $X$ is greater than $m$ and its homology
group in dimension $m$ is nontrivial, we derive that the homology
group of $\tilde{X}$ in dimension $m$ is also nontrivial. This
completes the proof of Proposition \ref{p:51}.
\end{pfn}

\section{Proof of Theorems \ref{t:14} and \ref{t:15}}
In this section we assume that assumptions $(A_0)$, $(A_5)$ and
$(A_6)$ hold and we are going to prove Theorems \ref{t:14} and
\ref{t:15}. First, we start by proving the following main results.
 \begin{pro}\label{p:61}
Let $n\geq 7$. There exists a pseudogradient $Y_2$ such that the
following holds:\\
 There exists a constant $c>0$ independent of $u=\sum_{i=1}^2\a_i
P\d_{(a_i,\l_i)}\in V(2,\e)$ such that
  $$(-\partial J(u),Y_2)_2\geq c\biggl(\e_{12}^{\frac{n-3}{n-4}}+\sum
\frac{1}{\l_i ^2}+\frac{\mid\n
K(a_i)\mid}{\l_i}+\frac{1}{(\l_id_i)^{n-3}}\biggr)\leqno{1)}$$
 $$(-\partial J(u+\ov{v}),Y_2+\frac{\partial \ov{v}}
 {\partial(\a_i,a_i,\l_i)}(Y_2))_2\geq c\biggl(\e_{12}^{\frac{n-3}{n-4}}+\sum
\frac{1}{\l_i ^2}+\frac{\mid\n
K(a_i)\mid}{\l_i}+\frac{1}{(\l_id_i)^{n-3}}\biggr)\leqno{2)}$$
 $3)$ $ Y_2$ is bounded and the only case where the maximum
of the $\l_i$'s increases along $Y_2$ is when the points $a_i$'s
are close to two different critical points $y_j$ and $y_r$ of $K$
with $-\D K(y_l)>0$ for $l=j,r$. Furthermore the least distance to
the boundary only increases if it is small enough.
 \end{pro}
\begin{pf}
We divide the set $V(2,\e)$ into three sets $A_1\cup A_2\cup A_3$
where, for $u=\sum\a_i P\d_{(a_i,\l_i)}\in V(2,\e)$,
$A_1=\{u/d_1\geq d_0 \mbox{ and }d_2\geq d_0\}$, $A_2=\{u/d_1\leq
d_0 \mbox{ and }d_2\geq 2d_0\}$, $A_3=\{u/d_1\leq 2d_0 \mbox{ and
}d_2\leq 2d_0\}$. We will build a vector field on each set and
then, $Y_2$ will be a convex combination of those vector fields.\\
{\bf 1st set } For $u\in A_1$. We can assume without loss of
generality that $\l_1\leq \l_2$. We introduce the following set
$T=\{i/\, \, |\n K(a_i)|\geq C_2/\l_i \}$ where $C_2$ is a large
constant.
The set $A_1$ will be divided into four subsets\\
 {\it 1st subset: } The set of $u$ such that
$\e_{12}\geq \frac{C_1}{\l_2^2} \, \mbox{ and } (10\l_1\geq \l_2
\, \mbox{ or } \mid\n K(a_1)\mid \geq \frac{C_2}{\l_1})$, where
$C_1$ is a large constant. In this case, we define $W_1$ as
 $$W_1=-M\l_2\frac{\partial P\d_2}{\partial
\l_2}+\sum_{i\in T}\frac{1}{\l_i}\frac{\partial P\d_i}{\partial
a_i} \frac{\n K(a_i)}{|\n K(a_i)|},$$
 where $M$ is a large
constant. Using Propositions \ref{p:33} and \ref{p:34}, we derive that
\begin{align}
(-\partial J(u),W_1)_2 & \geq M(c\e_{12}+O(\frac{1}{\l_2^2}))+
\sum_{i\in T}\biggl(\frac{|\n K(a_i)|}{\l_i}+O(\frac{1}{\l_i
^2}+\e_{12})\biggr)
\notag \\
 & \geq c\biggl(\e_{12}+\sum \frac{|\n K(a_i)|}{\l_i}+\frac{1}{\l_i
^2}\biggr).
\end{align}
{\it 2nd subset: } The set of $u$ such that $\e_{12}\geq
\frac{C_1}{\l_2^2} ,\, \,  10\l_1\leq \l_2 \mbox{ and } \mid\n
K(a_1)\mid \leq \frac{C_2}{\l_1}$. In this case, the point $a_1$
is close to a critical point $y$ of $K$. We define $W_2$ as
$$W_2=W_1+\sqrt{M}\l_1\frac{\partial P\d_1}{\partial\l_1}(sign(-\D
K(y))).$$ Using Propositions \ref{p:33} and \ref{p:34}, we obtain
\begin{align}
(-\partial J(u),W_2)_2 & \geq M(c\e_{12}+O(\frac{1}{\l_2^2}))
+\sqrt{M}(\frac{c}{\l_1^2}+O(\e_{12}))+ \sum_{i\in T}
(\frac{\mid\n K(a_i)\mid}{\l_i}+O(\frac{1}{\l_i
^2}+\e_{12})) \notag \\
 & \geq c(\e_{12}+\sum \frac{|\n K(a_i)|}{\l_i}+\frac{1}{\l_i
^2}).
\end{align}
{\it 3rd subset: } The set of $u$ such that $\e_{12}\leq
\frac{C_1}{\l_2^2}\,  \mbox{ and } (\mid\n K(a_1)\mid \geq
\frac{C_2}{\l_1}\, \mbox { or } \mid\n K(a_2)\mid \geq
\frac{C_2}{\l_2})$. In this case, the set $T$ is not empty, thus
we define
$$W_3'=\sum_{i\in T}\frac{1}{\l_i}\frac{\partial P\d_i}{\partial
 a_i} \frac{\n K(a_i)}{|\n K(a_i)|}.$$
 Using Proposition \ref{p:34}, we find
\begin{eqnarray}\label{e:43}
(-\partial J(u),W_3')_2  \geq c \sum_{i\in T}\biggl(\frac{|\n
K(a_i)|}{\l_i}+O(\frac{1}{\l_i ^2}+\e_{12})\biggr).
\end{eqnarray}
If we assume that ($\mid\n K(a_1)\mid\geq C_2/\l_1$ or $10\l_1\geq
\l_2$) and we choose $C_1<<C_2$, \eqref{e:43} implies the desired
estimate. In the other situation i.e. ($\mid\n K(a_1)\mid\leq
C_2/\l_1$ and $10\l_1\leq \l_2$), the point $a_1$ is close to a
critical point $y$ of $K$. As in the second case, we define
$W_3''$ as
$$W_3''=\frac{1}{\l_2}\frac{\partial P\d_2}{\partial a_2}
\frac{\n K(a_2)}{\mid\n K(a_2)\mid}+ \l_1\frac{\partial
P\d_1}{\partial\l_1}(sign(-\D K(y))).
$$
 Using Propositions \ref{p:33} and \ref{p:34}, we derive that
\begin{align}
(-\partial J(u),W_3'')_2 & \geq c (\frac{\mid\n
K(a_2)\mid}{\l_2}+O(\frac{1}{\l_2 ^2}+\e_{12}))+
c(\frac{1}{\l_1^2}+O(\e_{12}))\notag \\
 & \geq c(\e_{12}+\sum \frac{\mid\n K(a_i)
\mid}{\l_i}+\frac{1}{\l_i ^2}).
\end{align}
$W_3$ will be a convex combination of $W_3'$ and $W_3''$.\\
 {\it 4th subset: } The set of $u$ such that $\e_{12}\leq
\frac{C_1}{\l_2^2}\,  \mbox{ and } \mid\n K(a_i)\mid \leq
\frac{C_2}{\l_i} \mbox { for } i=1,2$. In this case, the
concentration points are near two critical points $y_i$ and $y_j$
of $K$. Two cases may occur: either $y_i=y_j$ or $y_i\ne y_j$.\\
- If $y_i=y_j=y$. Since $y$ is a nondegenerate critical point, we
derive that $\l_k\mid a_k-y\mid \leq c$ for $k=1,2$ and therefore
$\l_1\mid a_1-a_2\mid \leq c$. Thus we obtain $\e_{12}\geq
c(\l_1/\l_2)^{(n-4)/2}$ and therefore $\e_{12}\leq C_1/\l_2^2=
o(1/\l_1^2)$. In this case we define $W_4'=\l_1(\partial P\d_1/
\partial \l_1)(sign(-\D K(y)))$. Using Proposition \ref{p:33}, we
derive that
 \begin{eqnarray}
(-\partial J(u),W_4')_2\geq \frac{c}{\l_1^2}+O(\e_{12})\geq
c(\e_{12}+\sum \frac{\mid\n K(a_i)\mid}{\l_i}+\frac{1}{\l_i ^2}).
\end{eqnarray}
- If $y_i\ne y_j$. In this case we have $\e_{12}=o(1/\l_k^2)$ for
$k=1,2$. The vector field $W_4''$ will depend on the sign of $-\D
K(y_k)$, $k=i,j$. If $-\D K(y_i)<0$ ($y_i$ is near $a_1$), we
decrease  $\l_1$. If $-\D K(y_i)>0$ and $-\D K(y_j)<0$, we
decrease $\l_2$ in the case where $10\l_1\geq \l_2$ and we
increase $\l_1$ in the other case. If $-\D K(y_k)>0$ for $k=i,j$,
we increase both $\l_k$'s. Thus we obtain
 \begin{eqnarray}
 (-\partial J(u),W_4'')_2\geq c(\e_{12}+\sum
 \frac{\mid\n K(a_i)\mid}{\l_i}+\frac{1}{\l_i ^2}).
\end{eqnarray}
 The vector field $W_4$ will be a convex combination of $W_4'$ and
$W_4''$.\\
{\bf 2nd set } For $u\in A_2$, we have $|a_1-a_2|\geq d_0$.
Therefore $\e_{12}=o(1/\l_1)$ and $H(a_2,.)\leq c$. Let us define
$W_5=(1/\l_1)(\partial P\d_1/\partial a_1)(-\nu_1)$. Using
Proposition \ref{p:34}, we find
\begin{eqnarray}\label{e:47}
(-\partial J(u),W_5)_2\geq \frac{c}{\l_1}+O(\e_{12})
+\frac{c}{(\l_1d_1)^{n-3}}\geq \frac{c}{\l_1}
+\frac{c}{(\l_1d_1)^{n-3}}.
\end{eqnarray}
If $\l_1\leq 10\l_2$, then, in the lower bound of \eqref{e:47}, we
can make appear $1/\l_2$ and all the terms needed in 1). In the
other case i.e. $\l_1\geq 10\l_2$, we define $W_6$ as
$W_6=W_5+Y_1(P\d_2)$ and we obtain the desired estimate in this
case also.\\
{\bf 3rd set } For $u\in A_3$ i.e. $d_i\leq 2 d_0$ for $i=1,2$. We
have three cases.\\
{\it 1st case :} If there exists $i\in \{1,2\}$ (we denote by $j$
the other index) such that $M_1d_i\leq d_j$ where $M_1$ is a large
constant. In this case we define
\begin{eqnarray}
W_7=\sum \frac{1}{\l_i}\frac{\partial P\d_i}{\partial
a_i}(-\nu_i).
\end{eqnarray}
Using Proposition \ref{p:34}, we derive that
\begin{align}
\bigl(-\partial J(u),W_7\bigr) & \geq
c\sum_k\biggl(\frac{1}{\l_k}+\frac{1}{(\l_kd_k)^{n-3}}\biggr)+o(\e_{12}^\frac{n-3}{n-4})\\
 & +
O\biggl(\sum_k \frac{1}{\l_k}\bigg|\frac{\partial
\e_{12}}{\partial a_k} \bigg|
+\frac{1}{(\l_1\l_2)^{(n-4)/2}}\frac{1}{\l_k} \bigg|
\frac{\partial H(a_1,a_2)}{\partial a_k}\bigg|
+\l_k|a_1-a_2|\e_{12}^{\frac{n-1}{n-4}}\biggr)\notag
\end{align}
Since $M_1d_i\leq d_j$, then we have $|a_1-a_2|\geq d_j/2\geq
M_1d_i/2$. Thus we obtain
\begin{eqnarray}
\frac{1}{\l_k}\bigg|\frac{\partial \e_{12}}{\partial a_k}\bigg|+
+\frac{1}{(\l_1\l_2)^{(n-4)/2}}\frac{1}{\l_k}\bigg|\frac{\partial
H(a_1,a_2)}{\partial a_k}\bigg|+\e_{12}^{\frac{n-3}{n-4}}
=o\biggl(\sum_{r=1}^2\frac{1} {(\l_rd_r)^{n-3}}\biggr).
\end{eqnarray}
The same estimate holds for $\l_k |a_1-a_2|\e_{12}^{(n-1)/(n-3)}$.
Thus claim 1) follows in this case.\\
 {\it 2nd case :} If $d_2/M_1\leq d_1\leq M_1d_2$ and
$\l_2/M_2\leq \l_1\leq M_2\l_2$ where $M_2$ is chosen large
enough. In this case we define
\begin{eqnarray}
W_8=\frac{1}{\l_2}\sum_i\frac{\partial P\d_i}{\partial
a_i}(-\a_i\nu_i).
\end{eqnarray}
Using Proposition \ref{p:34} we derive that
\begin{align}\label{W8}
\bigl(-\partial J(u),W_8\bigr)_2 \geq &
\frac{c}{\l_2}\biggl(1+\sum_k\frac{1}{d_k(\l_kd_k)^{n-4}}+
c\a_1\a_2\frac{\partial \e_{12}}{\partial a_1}(\nu_1-\nu_2)\\
 & +\frac{c\a_1\a_2}{(\l_1\l_2)^{(n-4)/2}} \sum_k\frac{\partial
H(a_1,a_2)}{\partial a_k}\nu_k\biggr)
+o(\e_{12}^{\frac{n-3}{n-4}})\notag
\end{align}
Observe that $|\partial\e_{12}/\partial a_1||\nu_1-\nu_2|=
O(\e_{12})=o(1)$ and using the fact that $\partial
H(a_1,a_2)/\partial \nu_i\geq o((d_1d_2)^{(3-n)/2})$. It remains
to appear $\e_{12}$ in the lower bound. For this effect, if there
exists $i$ such that $\e_{12}\leq m/(\l_id_i)^{4-n}$ where $m$ is
a fixed large positive constant, then we can make appear $\e_{12}$
in \eqref{W8}. In the other case, we decrease both $\l_i$'s and we
define $W_9=-\sum\l_i\partial P\d_i/\partial \l_i$. Using
Proposition \ref{p:33}, we obtain
\begin{eqnarray}\label{W9}
\bigl(-\partial J(u),W_9\bigr)_2 \geq c\e_{12}+\sum_i
O\biggl(\frac{1}{\l_i ^2}+\frac{1}{(\l_id_i)^{n-4}}\biggr) \geq
c\e_{12}+\sum_i O\biggl(\frac{1}{\l_i ^2}\biggr).
\end{eqnarray}
Thus, in this case, we define the vector field as $W_8+W_9$. Using
\eqref{W8} and \eqref{W9}, we obtain the desired estimate.\\
 {\it 3rd case :} If $d_2/M_1\leq d_1\leq M_1d_2$ and there exists
$i$ (we denote $j$ the other index) such that $\l_i\geq M_2\l_j$.
In this case we increase $\l_j$, we decrease $\l_i$ and we move
the points along the inward normal vector. Then we define
$W_{10}=-2m\l_i\partial P\d_i/\partial \l_i+m\l_j\partial
P\d_j/\partial \l_j+W_7$ where $m$ is a large constant. Using
Propositions \ref{p:33} and \ref{p:34}, we derive that
\begin{align}
\bigl(-\partial J(u),W_{10}\bigr) \geq &
m\biggl(c\e_{12}+\frac{c}{(\l_jd_j)^{n-4}}
+O(\frac{1}{(\l_id_i)^{n-4}})\biggr)\\
 & +c\biggl(\sum \frac{1}{\l_k}
+ \frac{1}{(\l_kd_k)^{n-3}}+O(\e_{12})\biggr).\notag
\end{align}
Observe that, in this case, we have $\l_jd_j=o(\l_id_i)$ if we
choose $M_1/M_2$ so small. Thus the desired estimate follows.\\
 The proof of  Claim $1)$ is then completed. Claim $3)$ follows immediately from the
construction of $Y_2$. Claim $2)$ follows  from the estimate of $\ov{v}$
as in \cite{B1} and \cite{BCCH}.
\end{pf}\\
Now, arguing as in the proof of Proposition \ref{p:42}, we easily derive the following result.
\begin{cor}\label{c:62} Let $n\geq 7$.
The only critical points at infinity in $V(2,\e)$ correspond to
$P\d_{(y_i,\infty)}+P\d_{(y_j,\infty)}$ where $y_i$ and $y_j$ are
two different critical points of $K$ satisfying $-\D K(y_k)>0$ for
$k=i,j$. Such critical point has a Morse index equal to
$2n-\sum_{r=i,j} index(K,y_r)+1$.
\end{cor}
\begin{pro}\label{p:63} Let $n\geq 7$ and assume that $(P)$
has no solution. Then the following claims hold \\
i)\hskip 0.2cm \quad if $X=\ov{\cup_{y\in B}W_s(y)}$, where
$B=\{y\in \O/\n K(y)=0,\, -\D K(y)>0\}$, then $f_\l(C_{y_0}(X))$
retracts by deformation on $\cup_{y_i\in
X-\{y_0\}}W_u(y_0,y_i)_\infty\cup X_\infty$
 where $X_\infty=(\cup_{y_i\in
X}W_u(y_i)_\infty)$.\\
ii)\quad if $X=\ov{W_s(y_{i_0})}$, where $y_{i_0}$ satisfies
$$
K(y_{i_0})= \mbox{max }\{K(y_i)/\mbox{index }(K,y_i)=n-k, \quad
-\D K(y_i) >0 \}
$$
and if assumption $(A_7)$ holds, then $f_\l(C_{y_0}(X))$ retracts by deformation on \\
 $\cup_{y_i\in
X-\{y_0\}}W_u(y_0,y_i)_\infty\cup X_\infty \cup \sigma _1$
where $\sigma_1 \subset \cup_{y_i/index (K,y_i) \geq n-k}W_u(y_i)_\infty$.
\end{pro}
\begin{pf}
Let us start by proving Claim i).
Since $J$ does not have any critical point, the manifold
$f_\l(C_{y_0}(X))$ retracts by deformation on the union of the
unstable manifolds of the critical points at infinity dominated by
$f_\l(C_{y_0}(X))$ (see \cite{BR},\cite{M}). Proposition
\ref{p:42} and Corollary \ref{c:62} allow us to characterize such
critical points. Observe that we can modify the construction of
the pseudogradient defined in Proposition \ref{p:41} and
Proposition \ref{p:61} such that, when we move the point $x$ it
remains in $X$ i.e. we can use $Z_K$ instead of $\n K/\mid\n
K\mid$ where $Z_K$ is the pseudogradient for $K$ which we use to
build the manifold $X$. \\
For an initial condition $u=(\a/K(y_0)^{(n-4)/8})P\d_{(y_0,\l)}+
((1-\a)/K(x)^{(n-4)/8})P\d_{(x,\l)}$ in\\ $f_\l(C_{y_0}(X))$, the
action of the pseudogradient (see Proposition \ref{p:61}) is
essentially on $\a$. The action of bringing $\a$ to zero or to 1
depends on whether $\a<1/2$ (in this case, $u$ goes to $X_\infty$)
or $\a>1/2$ (in this case, $u$ goes to $\ov{W}_u((y_0)_\infty)$).
On the other hand we have another action on $x\in X$, when
$\a=1-\a=1/2$. Since only $x$ can move, then $y_0$ remains one of
the concentration points of $u$ and either $x$ goes to $W_s(y_j)$
where $y_j$ is a critical point of $K$ in $X-\{y_0\}$ or $x$ goes
to a neighborhood of $y_0$. In the last case  the flow has to exit from
$V(2,\e)$ (see the construction of
$Y_2$ in Proposition \ref{p:61}). The level of $J$ in this situation is close to
$(2S_n)^{4/n}/K(y_0)^{(n-4)/n}$ and therefore it cannot dominate
any critical point at infinity of two masses (since $K(y_0)=\max
K$). Thus the flow has to enter in $V(1,\e)$ and it will dominate
$(y_i)_\infty$ for $y_i\in X$.  Then $u$ goes to
$$
\biggl(\cup_{y_i\in
X-\{y_0\}}W_u((y_0,y_i)_\infty)\biggr)\cup\biggl(\cup_{y_i\in
X}W_u((y_i)_\infty)\biggr).
$$
 Then Claim i) follows. Now, using assumption $(A_7)$ and the same
 argument as in the proof of Claim i), we easily derive Claim ii).
 Thus our proposition follows.
\end{pf}\\
We now prove our theorems.\\
\begin{pfn}{\bf Theorem \ref{t:14} }
Arguing by contradiction, we assume that $(P)$ has no solution.
Using Proposition \ref{p:63} and the fact that  $\mu (y_{i_0})=0$,
we derive that $f_\l(C_{y_0}(X))$ retracts by deformation on
$X_\infty\cup D$ where $ D \subset \sigma$ is a stratified set  of
dimension at most $k$ ( in the topological sense, that is, $D\in
\Sig _j$, the group of chains of dimension $j$ with $j\leq k$) and
where $\sigma =\cup_{y_i\in X-(\{y_{i_0},y_0\})}
W_u((y_0,y_i)_\infty)\cup \cup_{y_i/index (K,y_i) \geq
n-k}W_u(y_i)_\infty$
is a manifold in dimension at most $k$. \\
As  $f_\l(C_{y_0}(X))$ is a contractible set, we then have
 $H_*(X_\infty\cup D) = 0$, for all $*\in \N^*$.
Using the exact homology sequence of $(X_\infty \cup D, X_\infty)$, we derive
$H_k(X_\infty)=H_{k+1}(X_\infty\cup D, X_\infty)=0$. This yields a contradiction
since $X_\infty \equiv X\times [A, +\infty)$, where $A$ is a large positive constant.
Therefore our theorem follows.
\end{pfn}\\
\begin{pfn}{\bf Theorem \ref{t:15}}
Assume that $(P)$ has no solution. By the above arguments, if
$\mu(y_i)=0$ for each $y_i\in B_k$, then $f_\l(C_{y_0}(X))$ retracts by
deformation  on $X_\infty\cup D$ where $ D\subset \sigma$ is a
stratified set and where  $\sigma =\cup_{y_i\in
X- (B_k\cup\{y_0\})} W_u((y_0,y_i)_\infty)$ is a manifold in
dimension at most $k$.\\
As in the proof of Theorem \ref{t:14}, we derive that
$H_*(X_\infty\cup D)=0$ for each $*$.  Using the exact homology sequence of
$(X_\infty\cup D, X_\infty)$ we obtain
$H_k(X_\infty)=H_{k+1} (X_\infty\cup D, X_\infty)=0$, this yields a
contradiction and therefore our result follows.
\end{pfn}
\small{

}
\end{document}